\documentclass[a4paper]{article}
\usepackage{a4wide}
\usepackage[utf8]{inputenc}
\usepackage[ngerman, english]{babel}
\usepackage{amsmath,amsthm,amssymb,amsfonts}
\usepackage{mathrsfs}
\usepackage{ifthen}
\usepackage{enumerate}
\usepackage{comment}
\usepackage{subcaption}
\usepackage{cleveref}
\usepackage{hyphsubst}
\usepackage{xcolor}
\usepackage{calc}
\usepackage{graphicx}
\usepackage[toc,page]{appendix}
\usepackage[bbgreekl]{mathbbol}
\DeclareSymbolFontAlphabet{\mathbbm}{bbold}
\DeclareSymbolFontAlphabet{\mathbb}{AMSb}
\usepackage{bbold}

\renewcommand{\Im}{\operatorname{Im}}

\newcommand{\indf}{\mathbbm{1}}

\makeatletter
\def\@hspace#1{\begingroup\setlength\dimen@{#1}\hskip\dimen@\endgroup}
\makeatother

\newtheorem{theorem}{Theorem}[section]
\newtheorem{lemma}[theorem]{Lemma}
\newtheorem{definition}[theorem]{Definition}

\newtheorem{remark}[theorem]{Remark}
\newtheorem{corollary}[theorem]{Corollary}

\newcommand{\dint}[1]{{#1}_{-}}
\newcommand{\dext}[1]{{#1}_{+}}
\newcommand{\trf}{u}
\newcommand{\tef}{v}
\newcommand{\rsf}{f}
\newcommand{\rsfa}{\tilde{f}}
\newcommand{\map}{\varphi}
\newcommand{\p}{x}
\newcommand{\q}{y}
\newcommand{\matB}{\boldB}
\newcommand{\matb}{\boldb}
\newcommand{\vecf}{\boldf}

\newcommand{\dd}{\mathrm{d}}

\newcommand{\ol}[1]{\overline{#1}}

\newcommand{\bpm}{\begin{pmatrix}}
\newcommand{\epm}{\end{pmatrix}}

\DeclareMathOperator{\dist}{dist}
\renewcommand{\div}{\operatorname{div}}

\DeclareMathOperator{\supp}{supp}

\newcommand{\setC}{\mathbb{C}}

\newcommand{\setN}{\mathbb{N}}

\newcommand{\setR}{\mathbb{R}}

\newcommand{\boldb}{\mathbf{b}}

\newcommand{\boldf}{\mathbf{f}}

\newcommand{\boldB}{\mathbf{B}}

\newcommand{\calC}{\mathcal{C}}
\newcommand{\calD}{\mathcal{D}}

\newcommand{\calL}{\mathcal{L}}

\newcommand{\calT}{\mathcal{T}}

\definecolor{brickred}{rgb}{0.8, 0.25, 0.33}
\definecolor{bostonuniversityred}{rgb}{0.8, 0.0, 0.0}
\definecolor{cornellred}{rgb}{0.7, 0.11, 0.11}
\definecolor{corn}{rgb}{0.98, 0.93, 0.36}
\definecolor{schoolbusyellow}{rgb}{1.0, 0.85, 0.0}
\definecolor{TUblue}{rgb}{0,102,153}
\colorlet{TUbluelight}{TUblue!30!white}

\usepackage{fancyhdr}
\usepackage[bbgreekl]{mathbbol}
\usepackage{authblk}
\title{A random walk in \LaTeX}
\author[1]{Martin Halla}
\author[1]{Thorsten Hohage}
\author[1]{Florian Oberender}
\affil{Institut f\"ur Numerische und Angewandte Mathematik,
Georg-August Universität Göttingen}
\title{A new numerical method for scalar eigenvalue problems in heterogeneous, dispersive, sign-changing materials\footnote{Support from DFG, CRC 1456 project 432680300 is gratefully acknowledged.}}

\begin{document}

\maketitle
\begin{abstract}
\noindent
We consider time-harmonic scalar transmission problems between dielectric and dispersive materials with generalized Lorentz frequency laws.
For certain frequency ranges such equations involve a sign-change in their principle part.
Due to the resulting loss of coercivity properties, the numerical simulation of such problems is demanding.
Furthermore, the related eigenvalue problems are nonlinear and give rise to additional challenges.
We present a new finite element method for both of these 
types of  problems, which is based on a weakly coercive reformulation of the PDE.
The new scheme can handle $C^{1,1}$-interfaces consisting piecewise of elementary geometries.
Neglecting quadrature errors, the method allows for a straightforward convergence analysis.
In our implementation we apply a simple, but nonstandard quadrature rule to achieve negligible quadrature errors.
We present computational experiments in 2D and 3D for both 
source and eigenvalue problems which confirm the stability and convergence of the new scheme.
\\

\\
\noindent%
\end{abstract}
\section{Introduction}\label{sec:introduction}
The starting point of this work are time-harmonic electromagnetic transmission problems involving dispersive materials, modeled by Maxwell's equations.
To simplify the setting we assume that the domain is invariant in one direction and bounded by a perfect conductor in the other two.
Thus the equations are reduced to two uncoupled systems called the transverse magnetic (TM) and transverse electric (TE) problem.
They can be further transformed into two scalar equations for the electromagnetic field (E,H) in the invariant direction with appropriate boundary conditions \cite{BonnetBDChesnelCiarlet:14a}.
Both equations have the form
\begin{align*}
	-\div(\sigma\nabla\trf)-\omega^{2}\tau \trf=\rsfa
\end{align*}
with the temporal frequency \(\omega\), the dispersive permeability $\mu=\mu(\omega,x)$ and permittivity $\epsilon=\epsilon(\omega,x)$, and \((\sigma,\tau)=(\mu^{-1},\epsilon)\) for the TE and \((\sigma,\tau)=(\epsilon^{-1},\mu)\) for the TM problem \cite{BonnetBDChesnelCiarlet:12}.
The case where \(\sigma\) is real valued and the sign of \(\sigma\) changes is of particular interest since then the arising bilinear forms are no longer (weakly) coercive, and classical theory fails \cite{BonnetBDChesnelCiarlet:14a}.
On the other hand, the term associated to 
$-\omega^{2}\tau(\omega,\cdot) \trf$ constitutes a compact perturbation for each $\tau(\omega,\cdot)\in L^\infty$.
Hence, this part is omitted for the Fredholmness and 
discretization analysis as its inclusion requires only 
standard arguments.

An important field where such sign-changing equations appear is the study of surface plasmons, which are electromagnetic waves that can form along the surface between a conductor and a dielectric material. They are the result of a resonance of light and free electrons on the surface of the conductor.  This resonance of electrons essentially traps the light along the surface and is only possible if the frequency-dependent permittivities of the two materials have different signs \cite{BDE03}.  These plasmons provide a unique way to concentrate and channel light and because of their properties there are many potential applications including light harvesting \cite{ALFSMP10}, the construction of miniaturized photonic circuits, and the detection of molecules \cite{BDE03}.


To obtain numerical solutions to this problem different strategies have been developed. In the case of piecewise constant coefficients, the problem can be solved using the boundary element method \cite{Unger:21}. 
Another approach suggested in 
\cite{AbdulleHuberLemaire:17} reduces the problem 
to a quadratic optimization problem for functions 
on the interface (see also the further development in \cite{CiarletLassounonRihani23} avoiding additional regularity assumptions). If the optimization problem is solved iteratively, e.g., by the conjugate gradient method, PDEs with coefficients of constant signs have 
to be solved in each iteration step. A further issue is the 
proper choice of a stabilization parameter.
Standard finite element methods in general only converge 
if the contrast (see \eqref{al:contrast} for a precise definition) is large enough as shown in  \cite{BonnetBDCiarletZwoelf:10} using T-coercivity techniques.
However, the necessary bounds for the contrast are not known explicitly and simulations computed this way can be treacherous as shown in \Cref{fig:ball_usual} and \Cref{fig:sol_wrong}.
Sharp convergence results, with respect to the contrast,
of finite element discretizations have been shown for polygonal  
interfaces in \cite{BonnetBDCarvalhoCiarlet:18} 
if angles are rational multiples of $2\pi$, meshes are choosen symmetric 
in a neighborhood of the flat parts, and special meshes at the 
corners are used.
While this approach is promising for 2D the construction of respective corner meshes in 3D is only possible for special cases (the Fichera corner) or leads to unsharp results.
The eigenvalue problems associated to dispersive transmission problems are nonlinear and have received significantly less attention so far.
See \cite{Unger:21} for BEMs and a \cite{Halla23SC} for a generalization of \cite{BonnetBDCarvalhoCiarlet:18} to 2D Maxwell EVPs.

In summary, the only method for elliptic differential equations 
with variable sign-changing coefficients for three-dimensional 
domains and non-polygonal two-dimensional domains is 
the optimization-based approach in \cite{AbdulleHuberLemaire:17}. However, the optimization process in  \cite{AbdulleHuberLemaire:17} requires a large number 
of PDE solutions, making this approach computationally 
significantly more costly than a finite element discretization,
especially when applied to eigenvalue problems.
The main aim of this work is to propose a finite element type 
discretization of such problems for smooth interfaces with only standard
requirements on the mesh. 
The principle approach of our new scheme is to apply a suitable T-operator to the PDE yielding a weakly coercive equation and allowing for discretization with standard 
finite element spaces. 
In addition, our method can naturally be applied to the 
related eigenvalue problems.

The remainder of this article is structured as follows.
First we specify the considered problem.
In \Cref{sec:reformulation} we introduce the applied reflection operator and the weakly coercive reformulation of the PDE.
In \Cref{sec:implementation} we discuss the implementation of the FEM and the used quadrature rules.
In \Cref{sec:numerical_exps} we present several computational experiments which confirm the stability and convergence of the new scheme.
In \Cref{sec:appendix} we include some technical analysis on the bounds of the used reflection operators.



\subsection*{Notations and problem setting}
For a domain \(U\subset\setR^{d}, d=2,3\) we denote by \((\cdot,\cdot)_{U}\) the scalar product of \(L^{2}(U)\) with associated norm \(\left\|\cdot\right\|_U\).
Furthermore we denote by \(\left\|\cdot\right\|_{H_{0}^{1}(U)}:=\left\|\nabla\cdot\right\|_{U}\) the norm of \(H_{0}^{1}(U)\) and by \(\left\|\cdot\right\|_{H^{-1}(U)}\) the norm on \(H^{-1}(U)\) that is given by the dual norm of \(H_{0}^{1}(U)\) with respect to the norm \(\left\|\nabla\cdot\right\|_{U}\). Unless specified otherwise, all spaces and scalar products are over \(\setR\). 

We consider a bounded Lipschitz domain \(\Omega\subset\setR^{d}\) which is decomposed into two disjoint, nonempty Lipschitz subdomains $\Omega_\pm\subset\Omega$ with a $C^{1,1}$-interface $\Gamma:=\partial\Omega_+\cap\partial\Omega_-$ such that $\ol{\Omega}_-\subset\Omega$.
We introduce the notation \(\dext{\trf}:=\trf|_{\dext{\Omega}}\) and \(\dint{\trf}:=\trf|_{\dint{\Omega}}\) for the restrictions of a function \(\trf\) defined on \(\Omega\).
Similarly, for a subdomain \(\Sigma\subset\Omega\) we write \(\dext{\Sigma}:=\Sigma\cap\dext{\Omega}\) and \(\dint{\Sigma}:=\Sigma\cap\dint{\Omega}\). We define the spaces of restrictions of \(H_{0}^{1}\) on these subdomains by
\begin{align*}
	H_{0,\Gamma}^{1}(\Sigma_\pm):=\{\trf|_{\Sigma_\pm}\colon\trf\in H_{0}^{1}(\Omega)\}.
\end{align*}
Additionally, we assume that the coefficient function \(\sigma\in L^{\infty}(\Omega)\), \(|\sigma|\) is essentially bounded from below by a positive constant and the restrictions of $\sigma$ satisfy
\(\dint{\sigma}>0\) and \(\dext{\sigma}>0\).
Lastly, we restrict ourselves to homogeneous Dirichlet boundary conditions and consider a source term \(\rsf\in H^{-1}(\Omega)\). This leads to the following problem:
\begin{align}\label{al:org_prob}
\begin{aligned}
	\text{Find } \trf \in H_{0}^{1}(\Omega) \text{ such that }
	-\div(\sigma\nabla \trf)=\rsf \text{ in }\Omega.
\end{aligned}
\end{align}
The former equation can alternatively be formulated in operator form using the operator \(B\in\calL(H_{0}^{1}(\Omega))\) defined by \((B\trf,\tef)_{H_{0}^{1}(\Omega)}:=(\sigma\nabla\trf,\nabla\tef)_{L^{2}(\Omega)}\):
\begin{align}\label{al:op}
\begin{aligned}
	\text{Find } \trf &\in H_{0}^{1}(\Omega) \text{ such that }
	B\trf=\rsfa
\end{aligned}
\end{align}
where \(\rsfa \in H_{0}^{1}(\Omega)\) is the image of \(\rsf\) under the canonical identification of \(H^{-1}(\Omega)\) and \(H_{0}^{1}(\Omega)\).

We also consider holomorphic eigenvalue problems related to the dispersive transmission problems.
We restrict ourselves to local lossless passive materials  which are nondispersive in $\Omega_+$.
Such materials are described by generalized Lorentz laws \cite[Theorem 3.22]{CassierJolyKachanovska17}, and convenient reconstructions from measurement data also take this form \cite{Grabovsky22}, cf.\ \cite{CassierJolyKachanovska17} for further discussions on the physical and mathematical requirements of dispersive material laws. 
Generalized Lorentz laws are described by coefficients $\sigma_0,\tau_0\in L^\infty(\Omega)$ such that $\sigma_0,\tau_0>0$ and $\sigma_0^{-1},\tau_0^{-1}\in L^\infty(\Omega)$, and take the form 
\begin{subequations}
\label{al:sigma_tau}
\begin{align}	\sigma(\omega,x)&=\left(\sigma_0(x)\bigg(1+\indf_{\Omega_-}(x)\sum_{l=1}^{N_\sigma}\frac{c_{\sigma,l}^2}{\omega_{\sigma,l}^2-\omega^2}\bigg)\right)^{-1}&\quad N_{\sigma}&\in\setN,&\omega_{\sigma,l},c_{\sigma,l}&\geq0,\\
	\tau(\omega,x)&=\tau_0(x)\bigg(1+\indf_{\Omega_-}(x)\sum_{l=1}^{N_\tau}\frac{c_{\tau,l}^2}{\omega_{\tau,l}^2-\omega^2}\bigg)&\quad N_{\tau}&\in\setN,&\omega_{\tau,l},c_{\tau,l}&\geq0.
\end{align}
\end{subequations}
Setting 
\begin{align*}
\Lambda:=\setC\setminus\bigg(
\bigcup_{l=1}^{N_\sigma}\{\pm\omega_{\sigma,l}\} \cup
\bigcup_{l=1}^{N_\tau}\{\pm\omega_{\tau,l}\} \cup
\{\omega\in\setC\colon\sigma(\omega)=0\}\bigg),
\end{align*}
we consider the following problem:
\begin{align}\label{al:eig}
\begin{aligned}
	\text{Find } (\omega,\trf)\in \Lambda\times H_{0}^{1}(\Omega,\setC)\setminus\{0\}& \text{ such that }
	-\div(\sigma(\omega)\nabla\trf)-\omega^2\tau(\omega)\trf=0 \text{ in }\Omega.
\end{aligned}
\end{align}
We rewrite the former as the eigenvalue problem for a holomorphic operator function:
\begin{align}\label{al:eig_op}
\begin{aligned}
    \text{Find } (\omega,\trf)\in \Lambda\times H_{0}^{1}(\Omega,\setC)\setminus\{0\}
	\text{ such that }
	B(\omega)\trf=0,
 \end{aligned}
 \end{align}
 with \(B(\cdot)\colon \Lambda\to\calL(H^1_0(\Omega,\setC))\) defined by
 \[(B(\omega)\trf,\tef)_{H_{0}^{1}(\Omega,\setC)}:=(\sigma(\omega,\cdot)\nabla\trf,\nabla\tef)_{L^{2}(\Omega,\setC^d)}
 -(\omega^2\tau(\omega,\cdot) \trf,\tef)_{L^2(\Omega,\setC)}.
 \]%
Note that for $\mp\Im(\omega^2)>0$ we have
 \[
	\pm\Im\left(\frac{-\omega^2}{\omega^2_{\tau,l}-\omega^2}\right)\geq0
	\qquad\text{and}\qquad
	\pm\Im(\sigma(\omega,\cdot))\geq0.
 \]
This shows that for a solution $(\omega,u)$ to \eqref{al:eig_op} with $\Im(\omega^2)\neq0$ it follows from $\Im (B(\omega)u,u)_{H_0^1(\Omega,\setC)}=0$ that $(\tau_0u,u)_{L^2(\Omega)}=0$ and hence $u=0$.
Since for $\Im(\omega^2)\neq0$ the operator $B(\omega)$ is weakly coercive, it follows from the Fredholm alternative that the spectrum of $B(\cdot)$ is real.
The challenging part is then to compute the part of the spectrum contained in $\Lambda_-:=\{\omega\in\setR\cap\Lambda\colon \sigma(\omega,\cdot)\not>0\}$.
Here $\sigma(\omega,\cdot)\not>0$ means that there exists a set 
$\Omega'\subset\Omega$ of positive measure such that 
$\sigma(\omega,x)\leq 0$ for all $x\in\Omega'$.

\section{The weakly coercive reformulation}\label{sec:reformulation}

As it is an important concept in our analysis, we start with a definition of (weak) coercivity of an operator.
\begin{definition}[(weak) coercivity]
An  operator \(B\in \calL(X)\) defined on a Hilbert space \(X\) is called coercive if there exists a constant \(\alpha>0\) such that
\[(B\tef,\tef)_{X}\geq\alpha\|\tef\|_{X}^{2}\quad\text{for all } \tef \in X.\]
An operator is called weakly coercive if it can be written as the sum of a compact and a coercive operator.
\end{definition}
It is well known that for weakly coercive operators or operator functions Galerkin schemes yield asymptotically reliable solutions  for source problems 
(see, e.g., \cite[(13.7b)]{K14}) or  for eigenvalue problems 
(\cite{Karma:96a,Karma:96b}), respectively.
For the problem at hand, the operator \(B\) is not weakly coercive. One technique to deal with problems lacking 
weak coercivity is the so called \(T\)-coercivity approach 
that we will briefly review now. The idea is to construct a bijective operator \(T\in\calL(H_{0}^{1}(\Omega))\) such that \(\tilde{B}:=T^{*}B\) is weakly coercive.
Then a solution to the source problem
\begin{align}\label{al:Tsource}
\begin{aligned}
	\text{Find } \trf &\in H_{0}^{1}(\Omega) \text{ such that }
	\tilde{B}\trf=T^{*}\rsfa
\end{aligned}
\end{align}
is a solution to the original problem \eqref{al:org_prob} and vice-versa.
For the eigenvalue problem we can proceed similarly if we are able to construct an operator \(T\) for which $T^*B(\omega)$ is weakly coercive for all $\omega\in\Lambda_-$.
In this case we consider the eigenvalue problem:
\begin{align}\label{al:Tevp}
\begin{aligned}
	\text{Find } (\omega,\trf) \in\Lambda'\times H_{0}^{1}(\Omega)\setminus\{0\} \text{ such that }
	T^{*}B(\omega)\trf=0,
\end{aligned}
\end{align}
where \(\Lambda'\subset\setC\) is an open neighborhood of \(\Lambda_-\) for which \(T^{*}B(\omega)\) is weakly coercive.
Since weak coercivity is a continuous property and \(T^{*}B(\omega)\) is weakly coercive on \(\Lambda_{-}\) such a neighborhood always exists, although to determine its exact shape an inspection of the frequency law is necessary.
Due to our assumptions \eqref{al:sigma_tau} the operator function \(B(\cdot)\) has only real eigenvalues and since \(T\) is bijective the problem above leads to the same eigenvalues and eigenfunctions as the original problem.

Since $T^*B$ or $T^*B(\omega)$ for $\omega\in\Lambda'$, respectively, are now weakly coercive, each Galerkin scheme yields an asymptotically converging approximation \cite[(13.7b)]{K14}, \cite{Karma:96a,Karma:96b}, and any
convenient finite element spaces can be used.

Several approaches have been suggested to construct \(T\)-operators yielding existence and convergence results, 
see \cite{BonnetBDCiarletZwoelf:10,BonnetBDChesnelCiarlet:12,BonnetBDCarvalhoCiarlet:18}.
However, the operators $T$ constructed in these references 
are not well suited for numerical implementations. 
Here we will work with a global reflection operator 
similar to \cite{BonnetBDCarvalhoCiarlet:18} for polygons in contrast to the patch-wise approach used in \cite{BonnetBDChesnelCiarlet:12}.
We define \(T\in\calL(H_{0}^{1}(\Omega))\) as either
\begin{align}\label{al:def_T}
	\dint{T}\trf:=\begin{cases}\dext{\trf}-2\chi\dint{R}\trf|_{\dint{\Sigma}}\quad&\text{in } \dext{\Omega}\\
	-\dint{\trf}\quad&\text{in } \dint{\Omega}\end{cases}
	\qquad\text{or}\qquad
	\dext{T}\trf:=\begin{cases}\dext{\trf}\quad&\text{in } \dext{\Omega}\\
	-\dint{\trf}+2\chi\dext{R}\trf|_{\dext{\Sigma}}\quad&\text{in } \dint{\Omega}\end{cases},
\end{align}
where \(\Sigma\) is a neighbourhood of \(\Gamma\),
\begin{align*}
	R_\pm\in\calL(H_{0,\Gamma}^{1}(\Sigma_\pm), H_{0,\Gamma}^{1}(\Sigma_\mp))
\end{align*}
are reflection operators which fulfill the so called matching condition \(\left(R_{\pm}u|_{\Sigma_{\pm}}\right)|_{\Gamma}=u_{\mp}|_{\Gamma}\) and \(\chi\in\calC^{1}(\Omega,[0,1])\) is a cut-off function with support in \(\Sigma\) which equals \(1\) in an open neighborhood of \(\Gamma\).
The weak coercivity of \(T^{*}B\) then depends on the operator norms of \(R_{\pm}\) and the so called contrasts of \(\sigma\) near the interface, that are given by
\begin{align}\label{al:contrast}
k_{+,\Sigma}:=\frac{\inf_{x\in\dext{\Sigma}}\sigma(x)}{\sup_{x\in\dint{\Sigma}}|\sigma(x)|} \quad\text{and}\quad
	k_{-,\Sigma}:=\frac{\inf_{x\in\dint{\Sigma}}|\sigma(x)|}{\sup_{x\in\dext{\Sigma}}\sigma(x)}.
\end{align}
Furthermore, we define \(k_{\pm}:=\inf_{\Sigma\supset\Gamma}k_{\pm,\Sigma}\) where the infimum is taken over all open neighborhoods of \(\Gamma\).
The precise relationship of the operators norms and the contrast is given by the following lemma,
the technique of which is well known (see, e.g., \cite{BonnetBDCiarletZwoelf:10,BonnetBDChesnelCiarlet:12}).
\begin{lemma}\label{lem:decomp}
For \(T_{\pm}\) be defined as above the following implication hold true:
\begin{align*}
\|R_{\pm}\|^{2}<k_{\pm,\Sigma}\quad\Rightarrow\quad T_{\pm}^{*}B\text{ is weakly coercive.}
\end{align*}
\end{lemma}
\begin{proof}
We will only prove the statement for \(\dint{T}\), and we will write \(T\) instead of \(\dint{T}\) for better readability. The statements for \(\dext{T}\) can be shown in the same way.
To show that \(T^{*}B\) is weakly coercive under the given assumption on the contrast, 
we define the operators \(A,K\in\calL(H_{0}^{1}(\Omega))\) via bilinear forms as
\begin{align*}
	(A\trf,\tef)_{H_{0}^{1}(\Omega)}&:=(|\sigma|\nabla\trf,\nabla\tef)_{\Omega}-2(\sigma\nabla\trf,\chi\nabla\dint{R}\tef|_{\dint{\Sigma}})_{\Omega}&\forall u,v&\in H_{0}^{1}(\Omega),\\
	(K\trf,\tef)_{H_{0}^{1}(\Omega)}&:=-2(\sigma\nabla\trf,\nabla\chi\dint{R}\tef|_{\dint{\Sigma}})_{\Omega}&\forall u,v&\in H_{0}^{1}(\Omega).
\end{align*}
Note that due to the boundedness of \(\dint{R}\) both bilinear forms are bounded and hence the operators are well-defined. Then \(A+K=T^{*}B\), and we will show that \(A\) is coercive and \(K\) is compact. For this, we write for brevity \(\dint{\tilde{R}}\tef\) instead of \(\dint{R}\tef|_{\dint{\Sigma}}\)
and use the subdivision of the domain to compute
\begin{align*}
	(A\trf,\trf)_{H_{0}^{1}(\Omega)}&=(|\sigma|\nabla\trf,\nabla\trf)_{\Omega\setminus\Sigma}+(\sigma\nabla\trf,\nabla \trf)_{\dext\Sigma}+(|\sigma|\nabla\trf,\nabla \trf)_{\dint{\Sigma}}-2(\sigma\nabla\trf,\chi\nabla\dint{\tilde{R}}\trf)_{\dext{\Sigma}}\\
	&\geq(|\sigma|\nabla\trf,\nabla\trf)_{\Omega\setminus\Sigma}+(|\sigma|\nabla\trf,\nabla \trf)_{\Sigma}-2(\sigma\nabla\trf,\chi\nabla\dint{\tilde{R}}\trf)_{\dext{\Sigma}}.
\end{align*}
In the same way as in the proof of Theorem 2.1 in \cite{BonnetBDChesnelCiarlet:12}, we choose \(\eta>0\) and apply Young's inequality to estimate
\begin{align*}
	(|\sigma|\nabla\trf,\nabla \trf)_{\Sigma}-2(\sigma\nabla\trf,\chi\nabla\dint{\tilde{R}}\trf)_{\dext{\Sigma}}\geq&(|\sigma|\nabla\trf,\nabla \trf)_{\Sigma}-\eta(\sigma\nabla\trf,\nabla\trf)_{\dext{\Sigma}}\\
	&-\frac{1}{\eta}(\sigma\chi\nabla\dint{\tilde{R}}\trf,\chi\nabla\dint{\tilde{R}}\trf)_{\dext{\Sigma}}\\
	\geq&(|\sigma|\nabla\trf,\nabla \trf)_{\Sigma}-\eta(\sigma\nabla\trf,\nabla\trf)_{\dext{\Sigma}}\\
	&-\frac{\big(\sup_{x\in\dext{\Sigma}}\sigma(x)\big)\|\dint{R}\|^2}{\eta}(\nabla\trf,\nabla\trf)_{\dint{\Sigma}}\\
	\geq&(1-\eta)\big(\inf_{x\in\dext{\Sigma}}\sigma(x)\big)(\nabla\trf,\nabla\trf)_{\dext{\Sigma}}\\
	&+\left(\inf_{x\in\dint{\Sigma}}|\sigma(x)|-\frac{\big(\sup_{x\in\dext{\Sigma}}\sigma(x)\big)\|\dint{R}\|^2}{\eta}\right)(\nabla\trf,\nabla\trf)_{\dint{\Sigma}}.
\end{align*}
Now we use \(\|\dint{R}\|^{2}<\frac{\inf_{x\in\dint{\Sigma}}|\sigma(x)|}{\sup_{x\in\dext{\Sigma}}\sigma(x)}\) and obtain 
\begin{align*}
	(|\sigma|\nabla\trf,\nabla \trf)_{\Sigma}-2(\sigma\nabla\trf,\nabla\dint{\tilde{R}}\trf)_{\dext{\Sigma}}&\geq\min((1-\eta)\inf_{x\in\dext{\Sigma}}\sigma,(1-\frac{1}{\eta})\inf_{x\in\dint{\Sigma}}|\sigma|)(\nabla\trf,\nabla\trf)_{\Sigma}.
\end{align*}
We can now fix \(\eta\in(0,1)\) sufficiently close to $1$ 
and find that there exists a constant \(\alpha_{1}>0\) such that
\[(|\sigma|\nabla\trf,\nabla \trf)_{\Sigma}-2(\sigma\nabla\trf,\nabla\dint{\tilde{R}}\trf)_{\dext{\Sigma}}\geq \alpha_{1}(\nabla\trf,\nabla \trf)_{\Sigma}.\]
This leads to
\begin{align*}
	(A\trf,\trf)_{H_{0}^{1}(\Omega)}&\geq(|\sigma|\nabla\trf,\nabla\trf)_{\Omega\setminus\Sigma}+\alpha_{1}(\nabla\trf,\nabla \trf)_{\Sigma}\\
	&\geq\min(\inf_{x\in\Omega\setminus\Sigma}|\sigma|,\alpha_{1})((\nabla\trf,\nabla\trf)_{\Omega\setminus\Sigma}+(\nabla\trf,\nabla \trf)_{\Sigma})=\alpha_{2}\|\trf\|_{H_{0}^{1}(\Omega)}^{2}
\end{align*}
with \(\alpha_{2}:=\min(\inf_{x\in\Omega\setminus\Sigma}|\sigma|,\alpha_{1})>0\) and shows that \(A\) is coercive.

For \(K\) we define the following operators to express \(K\) as a product of them. We write \(\iota\colon H_{0,\Gamma}^{1}(\dint{\Sigma})\rightarrow L^{2}(\dint{\Sigma})\) for the compact embedding operator and define \(M_{-2\sigma}\colon(L^{2}(\Omega))^{d}\rightarrow (L^{2}(\Omega))^{d}\) and \(M_{\nabla\chi}\colon L^{2}(\dint{\Sigma})\rightarrow (L^{2}(\Omega))^{d}\) as the multiplication operators with symbols \(-2\sigma\) and \(\nabla\chi\) respectively. Additionally, we write \(P_{\dint{\Sigma}}\colon H_{0}^{1}(\Omega)\rightarrow H_{0,\Gamma}^{1}(\dint{\Sigma})\) for the corresponding restriction operator. Because of the definitions of \(\sigma,\chi\) and \(H_{0,\Gamma}^{1}(\dint{\Sigma})\) all these operators are bounded, and we can now use the definition of \(K\) to obtain 
\begin{align*}
	(K\trf,\tef)_{H_{0}^{1}(\Omega)}&=(-2\sigma\nabla\trf,\nabla\chi\dint{R}\tef|_{\dint{\Sigma}})_{\Omega}\\
	&=(M_{-2\sigma}\nabla\trf,M_{\nabla\chi}\iota\dint{R}P_{\dint{\Sigma}}\tef)_{\Omega}\\
	&=(P_{\dint{\Sigma}}^{*}\dint{R}^{*}\iota^{*}M_{\nabla\chi}^{*}M_{-2\sigma}\nabla\trf,\tef)_{H_{0}^{1}(\Omega)}.
\end{align*}
This implies
\[K=P_{\dint{\Sigma}}^{*}\dint{R}^{*}\iota^{*}M_{\nabla\chi}^{*}M_{-2\sigma}\nabla,
\]
and because \(\iota\) is compact, so is \(\iota^{*}\). 
This means that \(K\) is compact because it is the product of compact and bounded operators. It follows that \(T^{*}B\) is weakly coercive. 
\end{proof}

Now we will provide an explicit construction of \(R_{\pm}\) based on the geometry of the interface and provide upper bounds for their norm to clarify how they have to be constructed to achieve weak coercivity.

\subsection{Global reflection operators}\label{sec:global_reflection}
We construct \(R_\pm\) via a $C^{0,1}$-homomorphism
\(\map\colon\Sigma\rightarrow\Sigma\) with \(\map(\Sigma_{\pm})=\Sigma_{\mp}\) 
for some neighborhood $\Sigma$ of $\Gamma$ such that 
$\varphi(x)=x$ for all $x\in\Gamma$ by
\begin{align*}
	\dint{R}w&:=w\circ\map,\qquad
	\dext{R}w:=w\circ\map.
\end{align*}
Recall that if $\varphi\in C^{0,1}(\Sigma)$, then by the 
Rademacher theorem the Jacobian $D\varphi(x)$ exists for 
almost all $x\in \Sigma$, and $|\varphi|_{C^{0,1}}:=\mathrm{ess\,sup}_{x\in\Sigma}|D\varphi(x)|_{\mathcal{L}(\setR^2)}<\infty$. 
Moreover, if $w\in H^1(\Sigma)$, then $w\circ \varphi\in H^1(\Sigma)$ (see \cite[Theorem 4.1]{wloka:87}) 
and 
\begin{align}\label{eq:composition_bound}
\|w\circ \varphi\|_{H^1}\leq |\varphi|_{C^{0,1}}\|w\|_{H^1}.
\end{align}
To explicitly construct \(\map\) under our assumptions, we can define the unit normal vector \(n\colon\Gamma\rightarrow S^{d-1}\) pointing towards \(\dint{\Omega}\) everywhere on the interface and consider the functions
\begin{align}\label{eq:defiPhiM}
\begin{aligned}
	\Phi&\colon\Gamma\times(-\delta,\delta)\rightarrow\Omega,
	&&(\p,t)\mapsto \p+tn(\p),\\
	M&\colon\Gamma\times\setR\rightarrow\Gamma\times\setR,
	&&(\p,t)\mapsto(\p,-t).
\end{aligned}
\end{align}
Since $\Gamma$ is $C^{1,1}$-smooth by assumption, it 
follows that $n\in C^{0,1}(\Gamma,\setR^d)$, and hence 
$\Phi\in C^{0,1}(\Sigma,\setR^d)$.
If we can choose \(\delta>0\) small enough such that \(\Phi\) is a $C^{0,1}$-homomorphism, then we can define \(\map\) by
\[
	\map:=\Phi\circ M\circ\Phi^{-1}.
\]
Subsequently, we define \(\Sigma:=\Phi(\Gamma\times(-\delta,\delta))\) and note that \(\map(\Sigma_{\pm})=\Sigma_{\mp}\).
As \(M=M^{-1}\), we have  \(\map=\map^{-1}\).
For many 
practically relevant surfaces such as arcs, lines, planes and parts of spheres and cylinders an explicit computation of 
\(\Phi^{-1}\) is feasible, but for general surfaces, one has 
to resort to numerical inversions.
It is now possible to calculate upper bounds for the norms of \(\dext{R}\) and \(\dint{R}\) which only depend on the geometry of the interface and \(\delta\).

\begin{theorem}\label{the:bounds}
In two dimensions the reflection operators are bounded in norm by
\begin{align}
	\|\dint{R}\|&\leq\max\left(1,\left|\frac{1-\delta\inf_{\p\in\Gamma}\kappa(\p)}{1+\delta\inf_{\p\in\Gamma}\kappa(\p)}\right|\right),\qquad
	\|\dext{R}\|\leq\max\left(1,\left|\frac{1+\delta\sup_{\p\in\Gamma}\kappa(\p)}{1-\delta\sup_{\p\in\Gamma}\kappa(\p)}\right|\right),
\end{align}
where \(\kappa\) is the curvature of the interface.
In three dimensions the bounds are
\begin{align*}
	\|\dint{R}\|&\leq\max\left(1,\left|\frac{1-\delta\inf_{\p\in\Gamma}\kappa_{1}(\p)}{1+\delta\inf_{\p\in\Gamma}\kappa_{1}(\p)}\right|,\left|\frac{1-\delta\inf_{\p\in\Gamma}\kappa_{2}(\p)}{1+\delta\inf_{\p\in\Gamma}\kappa_{2}(\p)}\right|\right),\\
	\|\dext{R}\|&\leq\max\left(1,\left|\frac{1+\delta\sup_{\p\in\Gamma}\kappa_{1}(\p)}{1-\delta\sup_{\p\in\Gamma}\kappa_{1}(\p)}\right|,\left|\frac{1+\delta\sup_{\p\in\Gamma}\kappa_{2}(\p)}{1-\delta\sup_{\p\in\Gamma}\kappa_{2}(\p)}\right|\right)
\end{align*}
where \(\kappa_{1},\kappa_{2}\) are the principal curvatures of the interface.
\end{theorem}
The proof of these bounds relies on an application of the transformation formula and basic differential geometry.
It is given in Appendix \ref{sec:appendix}.
Using these explicit bounds, we can formulate conditions on the contrast of \(\sigma\) under which the construction of an operator \(T\) based on global reflection operators is possible such that \(T^{*}B\) is weakly coercive. For this we just combine our previous results with the fact that the bounds for the reflection operators decay to 1 when \(\delta\) gets smaller.
\begin{theorem}[Conditions for weak coercivity]\label{the:conditions}
For an interface that is \(C^{1,1}\) and piecewise \(C^{2}\) there exists an operator \(T\) which can be constructed via a global reflection operator as in \eqref{al:def_T} such that \(T^{*}B\) is weakly coercive if the following conditions are satisfied:
\begin{enumerate}
	\item There exists \(\delta_{0}>0\) such that the map \(\Phi\colon\Gamma\times(-\delta_{0},\delta_{0})\rightarrow\Omega\) defined by \eqref{eq:defiPhiM} 
	is a $C^{0,1}$-homomorphism. 
\item The curvature of \(\Gamma\) for the two dimensional case or the two principal curvatures of \(\Gamma\) for the three dimensional case are bounded.
\item One of the contrasts \(\dext{k}\) or \(\dint{k}\) of \(\sigma\) is strictly greater than \(1\).
\end{enumerate}
\end{theorem}
\begin{proof}
We only prove the case where \(\dext{k}>1\). The other case can be proven in the same way. From the second condition we know that the principal curvatures or the curvature is bounded in absolute value by a constant which we call \(\kappa_{max}\).
Now, since \(\dext{k}>1\) we can choose \(\delta\in(0,\delta_{0})\) small enough such that
\[
	1\leq\frac{1+\kappa_{max}\delta}{1-\kappa_{max}\delta}<\sqrt{\dext{k}}.
\]
Next we consider the map \(\Phi_{\delta}\colon\Gamma\times(-\delta,\delta)\rightarrow\Omega,(\p,t)\mapsto \p+tn(\p)\) which is a $C^{0,1}$-homomorphism due to the first assumption. 
Therefore, as in \Cref{sec:global_reflection}, we can  
construct the global reflection operator \(\dext{R}\). Then we can use \Cref{the:bounds} and obtain 
\[\|\dext{R}\|\leq\max\left(1,\left|\frac{1+\delta\sup_{\p\in\Gamma}\kappa(\p)}{1-\delta\sup_{\p\in\Gamma}\kappa(\p)}\right|\right)\leq\frac{1+\kappa_{max}\delta}{1-\kappa_{max}\delta}<\sqrt{\dext{k}}\]
for the two dimensional case and
\[\|\dext{R}\|\leq\max\left(1,\left|\frac{1+\delta\sup_{\p\in\Gamma}\kappa_{1}(\p)}{1-\delta\sup_{\p\in\Gamma}\kappa_{1}(\p)}\right|,\left|\frac{1+\delta\sup_{\p\in\Gamma}\kappa_{2}(\p)}{1-\delta\sup_{\p\in\Gamma}\kappa_{2}(\p)}\right|\right)\leq\frac{1+\kappa_{max}\delta}{1-\kappa_{max}\delta}<\sqrt{\dext{k}}\]
for the three dimensional case. Now \Cref{lem:decomp} implies that the operator \(\dext{T}^{*}B\) is weakly coercive with \(\dext{T}\) constructed via \(\dext{R}\). 
\end{proof}
\begin{remark}
From the proof we can see that the last condition can be slightly weakened. In general it is enough to require one of the contrasts \(k_{+,\Sigma}\) or \(k_{-,\Sigma}\) to be bounded from below by \(1\) where \(\Sigma\) contains all points that are closer to \(\Gamma\) than a fixed distance \(\delta>0\) which can be arbitrary small. Furthermore, it may seem to be advantageous to choose \(\delta\) as small as possible, but this comes with the price of a large gradient of the cut-off function.
\end{remark}
\begin{remark}
We also note that the required bounds on the contrast in the two dimensional case coincides with the ones 
in \cite{BonnetBDChesnelCiarlet:12} where the optimality 
of these bounds has been shown. 
\end{remark}
For simple geometries where the interface consists of circular arches and straight lines in two dimensions or planes, parts of spheres and cylinders in three dimensions the operator \(T\) can be implemented and used for finite element methods. Precise bounds for the necessary size of \(\delta\) for a given contrast are presented in Appendix \ref{sec:appendix}. Using such a suitable \(\delta\) the convergence of the source and eigenvalue problem is then established by the weak coercivity.

However there is one further challenge. For the full discretization the entries of the system matrix have to be computed numerically where integrals are approximated by quadrature rules. Usually this does not pose a major problem as long as the finite element functions and coefficient functions are smooth enough, because the Bramble-Hilbert lemma can be used to show that the quadrature error converges to zero for decreasing mesh sizes. Unfortunately, this is not the case for this method, because here we also have to consider integrals of finite element basis functions which have non-intersecting supports and numerically approximate integrals of the form
\[\int_{\calD}\sigma (\nabla u)^\top \nabla(v\circ\map)\dd x
=\int_{\calD}\sigma (\nabla u)^\top \nabla(v\circ\map)\indf_{\map(\supp v)}\dd x,\]
where we have left out the cut-off function \(\chi\) for simplicity. The problem with the numerical approximation of such integrals is that even in the simplest case where \(\map\) is an affine transformation and \(\nabla u\) and \(\nabla v\) are constant, the function \(\indf_{\map(\supp v)}\) is only in \(L^{\infty}(\calD)\) and therefore the classical methods fail. Additionally even if the jump of this function only occurs along a polygonal line, the quadrature approximation still does not get better with decreasing \(h\) because the function gets scaled as well. 
We therefore cannot hope to achieve convergence of the quadrature error to zero for decreasing \(h\). 
This will be further discussed in the forthcoming
section.

\section{Implementation}\label{sec:implementation}
We implemented our method using the finite element library NGSolve \cite{Schoeberl:14}.
The main effort lies in the implementation of the custom assemble procedure for the calculation of the stiffness matrix.
To explain its details we recall in the following the convenient framework of a finite element implementation.

Let \(\calT_{h}=\bigcup_{m}\calD_{m}\) be a mesh consisting of elements \(\calD_{m}\). Even though it is not necessary for our implementation, we assume that all elements are images of a  single reference element \(\calD\) to simplify the presentation. The FEM code then provides us with transformations \(\Psi_{m}:\calD\rightarrow\calD_{m}\) for each element as well as their Jacobians. The finite element space \(V_{h}\) is then implemented by providing a collection of shape functions \(s_{j}:\calD\rightarrow\setR\) for \(1\leq j \leq N\)  on the reference element.
For \(H^{1}\)-finite elements we can additionally access the gradients of the shape functions to calculate the gradient of a finite element function via the chain rule. Finally the FEM code also provides quadrature points \(\q_{i}\in\calD\) for \(1\leq j \leq N\) on the reference element.
With this tools we are able to outline the calculation of the stiffness matrix \(\matB:=(\matb_{i,j})_{i,j=1,\dots,P}\) and the right side vector \(\vecf=(\vecf_i)_{i=1,\dots,P}\) defined by
\begin{align*}
	\matb_{i,j}&:=(B\tef_{j},T\tef_{i})_{H^1_0(\Omega)}=\int_{\Omega}\sigma(\nabla\tef_{j})^{\top}\nabla(T\tef_{i})\dd x\\
	\vecf_{i}&:=\int_{\Omega}\rsfa \,T\tef_{i}\dd x
\end{align*}
for finite element functions \(\tef_{i},\tef_{j} \in V_{h}\subset H_{0}^{1}(\Omega)\). In the following part we will only consider the case where \(T\) is defined by
\begin{align*}
	\dint{T}\trf:=\begin{cases}\dext{\trf}-2\chi\dint{R}\trf|_{\dint{\Sigma}}\quad&\text{on } \dext{\Omega}\\
	-\dint{\trf}\quad&\text{on } \dint{\Omega}\end{cases},
\end{align*}
because the implementation of the other case is essentially the same. We will also write \(R\) instead of \(\dint{R}\). With this definition of \(T\) we can write
\begin{align*}
	\matb_{i,j}&=\int_{\Omega}\sigma(\nabla\tef_{j})^{\top}\nabla(T\tef_{i})\dd x\\
	&=\int_{\dext{\Omega}}\sigma(\nabla\tef_{j})^{\top}\nabla(\tef_{i}-2\chi R\tef_{i}|_{\dint{\Sigma}})\dd x+\int_{\dint{\Omega}}\sigma(\nabla\tef_{j})^{\top}\nabla(-\tef_{i})\dd x\\
	&=\int_{\Omega}|\sigma|(\nabla\tef_{j})^{\top}\nabla\tef_{i}\dd x-2\int_{\dext{\Sigma}}\sigma(\nabla\tef_{j})^{\top}\nabla(\chi R\tef_{i}|_{\dint{\Sigma}})\dd x\\
	&=:\matb^{(1)}_{i,j}-2\matb^{(2)}_{i,j}
\end{align*}
so we get \(\matB=\matB^{(1)}-2\matB^{(2)}\). In the same way we can define \(\vecf^{(1)},\vecf^{(2)}\) with \(\vecf=\vecf^{(1)}-2\vecf^{(2)}\). We note that \(\matB^{(1)}\) is the stiffness matrix of a bilinear form without any special operators, so it can be calculated the usual way.
The same is true for \(\vecf^{(1)}\). Additionally, we see that the domain of integration required for \(\matB^{(2)}\) is just \(\dext{\Sigma}\), so for its calculation we only have to consider finite element functions with support in \(\Sigma\). To make use of this and to generally simplify further calculations, we subdivide \(\Sigma=:\bigcup_{l=1}^{N}\Sigma^{(l)}\) according to the type of the interface by lines or planes which are perpendicular to the interface. We then generate the mesh such that each element lies entirely in either \(\Sigma_{l}\cap\dint{\Sigma}\) or \(\Sigma_{l}\cap\dext{\Sigma}\). 
This subdivision allows us to write \(R\) as
\begin{align*}
	(R\tef)(x)=(\tef\circ\map^{(l)})(x) \text{ for } x\in\dext{\Sigma}^{(l)}
\end{align*}
where \(\map^{(l)}\colon\Sigma^{(l)}\rightarrow\Sigma^{(l)}\) is a predefined transformation based on the interface geometry. Because the mesh respects this subdivision, on each mesh cell only one transformation has to be considered and transformations only have to be considered for mesh cells in \(\Sigma\).

Due to the presence of the transformation \(\map\) the assemble procedure of \(\matB^{(2)}\) is non-standard and poses the main challenge for the application of the method. This is mainly caused by the fact, that usual assemble procedures assume that only finite element functions which are supported on a common mesh cell contribute to the entries of the matrix. In our case functions \(\tef\) and \(\trf\circ\map\) can have intersecting supports even though the corresponding finite element functions have not. Additionally the intersection of the supports is usually not a mesh cell. To tackle these problems 
the main idea is to reduce the problem to contributions from single quadrature points which can then be calculated explicitly. For this let \(\calD_{m}\) be an element of \(\calT_{h}\) in \(\dext{\Sigma}\) and let \(\trf,\tef\in V_{h}\) be such that \(\supp\trf\cap\supp\tef\circ\map\cap\calD_{m}\neq\emptyset\). We can then assume that \(\trf\) is given on \(\calD_{m}\) by a single shape function and can therefore be written as \(\trf|_{\calD_{m}}=s^{(u)}\circ\Psi_{m}^{-1}\).
The contribution to the entry \(\matb_{(\tef,\trf)}\) corresponding to \(\trf\) and \(\tef\) is then given by
\begin{align*}
    \matb_{(\tef,\trf,m)}&:=\int_{\calD_{m}}\sigma(\p)\nabla(s^{(\trf)}\circ\Psi_{m}^{-1})(\p)^{\top}\nabla(\chi\cdot R\tef)(\p)\dd \p\\
    &=\int_{\calD}\sigma(\Psi_{m}(\q))\nabla(s^{(\trf)}\circ\Psi_{m}^{-1})(\Psi_{m}(\q))^{\top}\nabla(\chi\cdot R\tef)(\Psi_{m}(\q))|\det D_{\q}\Psi_{m}|\dd \q,
\end{align*}
where we write \(D_{\q}\Psi\) for the Jacobian of \(\Psi_{m}\) at \(\q\).
Because we are now integrating over the reference element, we can use the quadrature points to approximate the integral. If we now consider a single quadrature point \(\tilde{\q}\) with corresponding weight \(\tilde{w}\) then the contribution of this point is
\begin{align*}
	\matb_{(\tef,\trf,m,\tilde{\q})}&:=\tilde{w}\sigma(\Psi_{m}(\tilde{\q}))\nabla(s^{(\trf)}\circ\Psi_{m}^{-1})(\Psi_{m}(\tilde{\q}))^{\top}\nabla(\chi\cdot R\tef)(\Psi_{m}(\tilde{\q}))|\det D_{\tilde{\q}}\Psi_{m}|.
\end{align*}
For this single point we can find the element to which \(\Psi_{m}(\tilde{\q})\) is mapped by \(\map\) and call it \(\calD_{n}\). On this element \(\tef\) can then be represented by a single shape function \(s^{(\tef)}\). To shorten the expressions a little bit we combine scalar terms and define \(C_{m,\tilde{\q}}:=\tilde{w}|\det D_{\tilde{\q}}\Psi_{m}|\). Now we can use the definition of \(R\) to get
\begin{align*}
    \matb_{(\tef,\trf,m,\tilde{\q})}&=C_{m,\tilde{\q}}\sigma(\Psi_{m}(\tilde{\q}))\left(\nabla(s^{(\trf)}\circ \Psi_{m}^{-1})(\Psi_{m}(\tilde{\q}))\right)^{\top}\nabla(\chi\cdot (\tef\circ\map))(\Psi_{m}(\tilde{\q}))\\
    &=C_{m,\tilde{\q}}\sigma(\Psi_{m}(\tilde{\q}))\left(\nabla(s^{(\trf)}\circ \Psi_{m}^{-1})(\Psi_{m}(\tilde{\q}))\right)^{\top}\nabla\left(\chi\cdot(s^{(\tef)}\circ \Psi_{n}^{-1}\circ\map)\right)(\Psi_{m}(\tilde{\q}))
\end{align*}
which can be explicitly calculated via the chain rule as long as the values and derivatives of \(\map\) and \(\chi\) can be computed. The assembly of the whole matrix \(\matB^{(2)}\) then just consist of the calculation and summation of these terms together with some bookkeeping to add the contributions at the right position. Additionally all shape functions on an element can be considered at the same time and the implementation is done in a way that the calculation of \(\matb_{(\tef,\trf,m,\tilde{\q})}\) can be easily replaced to allow for the computation of different integrands. For example 
\[\matb_{(\tef,\trf,m,\tilde{\q})}=C_{m,\tilde{\q}}\left((s^{(\trf)}\circ \Psi_{m}^{-1})(\Psi_{m}(\tilde{\q}))\right)\cdot\left(\chi\cdot(s^{(\tef)}\circ \Psi_{n}^{-1}\circ\map)\right)(\Psi_{m}(\tilde{\q}))\]
can be used for \(L^{2}\)-terms.

\begin{figure}
\centering
       \includegraphics[clip,scale=0.1]{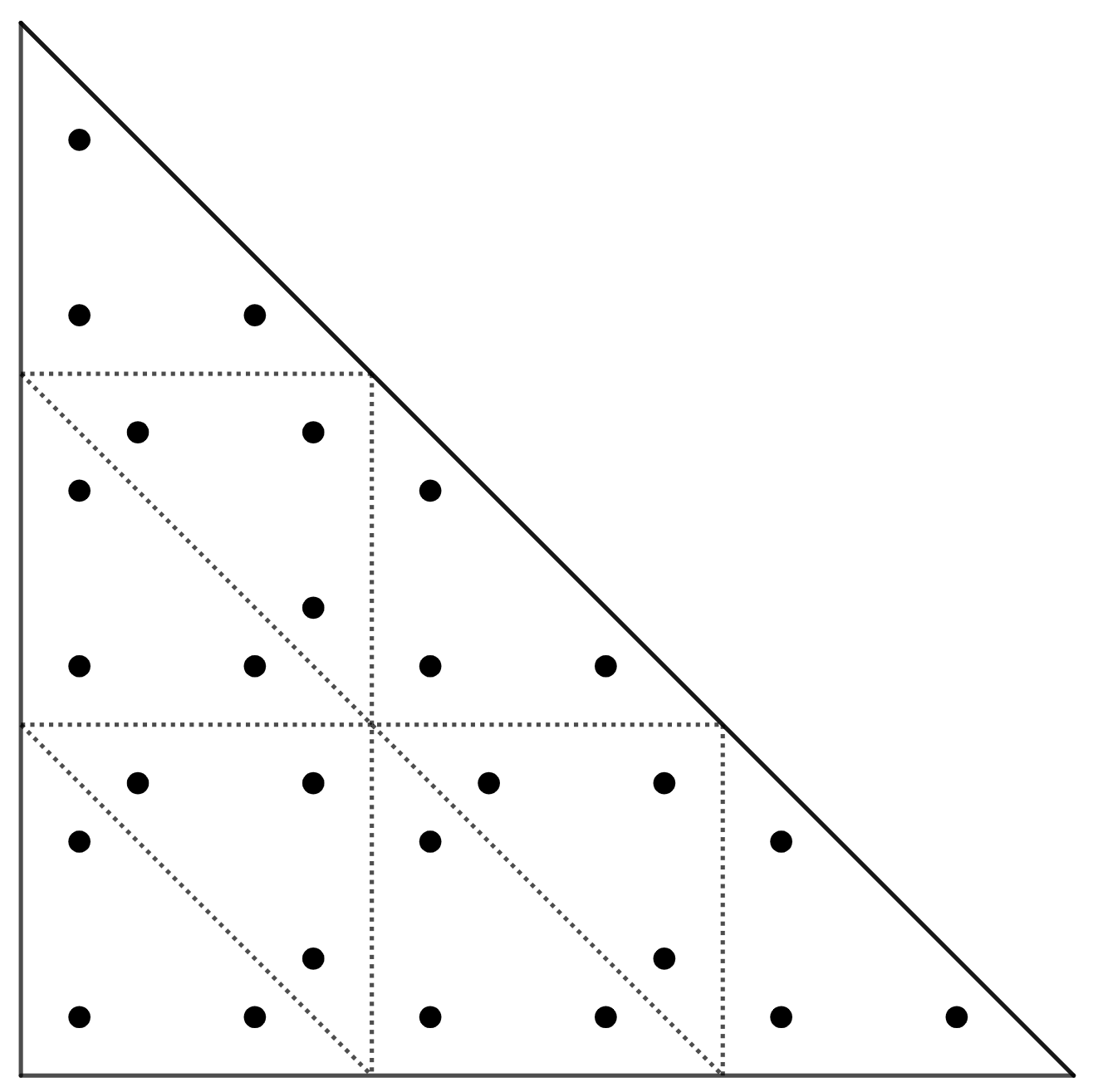}
\caption{Examples for subdivided element with quadrature points}
\label{fig:quadrature}
\end{figure}

To cope with the difficulties relating to the quadrature one needs to use a high number of quadrature points. This could be done by just using higher quadrature orders but they are adapted to high order polynomials and do not work particularly well for piece-wise continuous functions. We therefore use a hybrid quadrature rule which is based on dividing an element into many smaller similar elements and then using a standard Gauss-Legendre quadrature rule on each smaller element. An example for this is shown in \Cref{fig:quadrature}.
For the linear forms we can use a simple transformation to significantly simplify the calculations: 
\begin{align*}
	\boldf_{i}^{(2)}&=\int_{\dext{\Sigma}}\rsfa(\p)\chi(\p)(\tef_{i}\circ\map)(\p)\dd \p\\
	&=\int_{\dint{\Sigma}}\rsfa(\map(\q))\chi(\map(\q))\tef_{i}(\q)|\det D_{\q}\map|\dd \q.
\end{align*}
In this new form the transformation \(\map\) is now no longer composed with a finite element function so the calculation of \(\boldf^{(2)}\) can then be performed by implementing a coefficient function \(h(\q):=\rsfa(\map(\q))\chi(\map(\q))|\det D_{\q}\map|\) and then defining and assembling the linear form as usual.

\section{Numerical experiments}\label{sec:numerical_exps}
In this section we present different examples that illustrate our method for different domains. We consider the convergence rates, analyze the errors, and
for settings which allow the application of classical finite elements we compare our results with those.

\subsection{Examples in two dimensions}
Our first example is a two dimensional domain, that consists of a ring and a disc by defining \(\sigma:B_{2}\rightarrow\setR\) in polar coordinates by 
\begin{align*}
	\sigma(r):=\begin{cases}-1&\text{ for }r\leq1\\3&\text{ else}\end{cases},
\end{align*}
where $B_r:=\{x\in\setR^2\colon|x|<r\}$.
This leads to the subdomains \(\dint{\Omega}=B_{1}\subset\setR^{2}\) and \(\dext{\Omega}=B_{2}\setminus\overline{B_{1}}\).
We define the right hand-side \(\rsf\) also in polar coordinates as
\[\rsf(r):=\sigma(r)(9r-6),\]
which corresponds to the solution \(\trf_{\text{ref}}\in H_{0}^{1}(\Omega)\) given by
\begin{align}\label{eq:uref2D}
\trf_{\text{ref}}(r)=r^{3}-\frac{3}{2}r^{2}-2.
\end{align}
According to \Cref{cor:curved} we can choose \(\delta=0.2\) 
in \eqref{eq:defiPhiM} because the squared norm of \(R\) is then bounded by $2.25$, which is smaller than the contrast \(k_{-}=3\). 

\begin{figure}
     \centering
     \begin{subfigure}[b]{0.45\textwidth}
         \centering
        \includegraphics[trim={0 0 0 0},clip,scale=0.4]{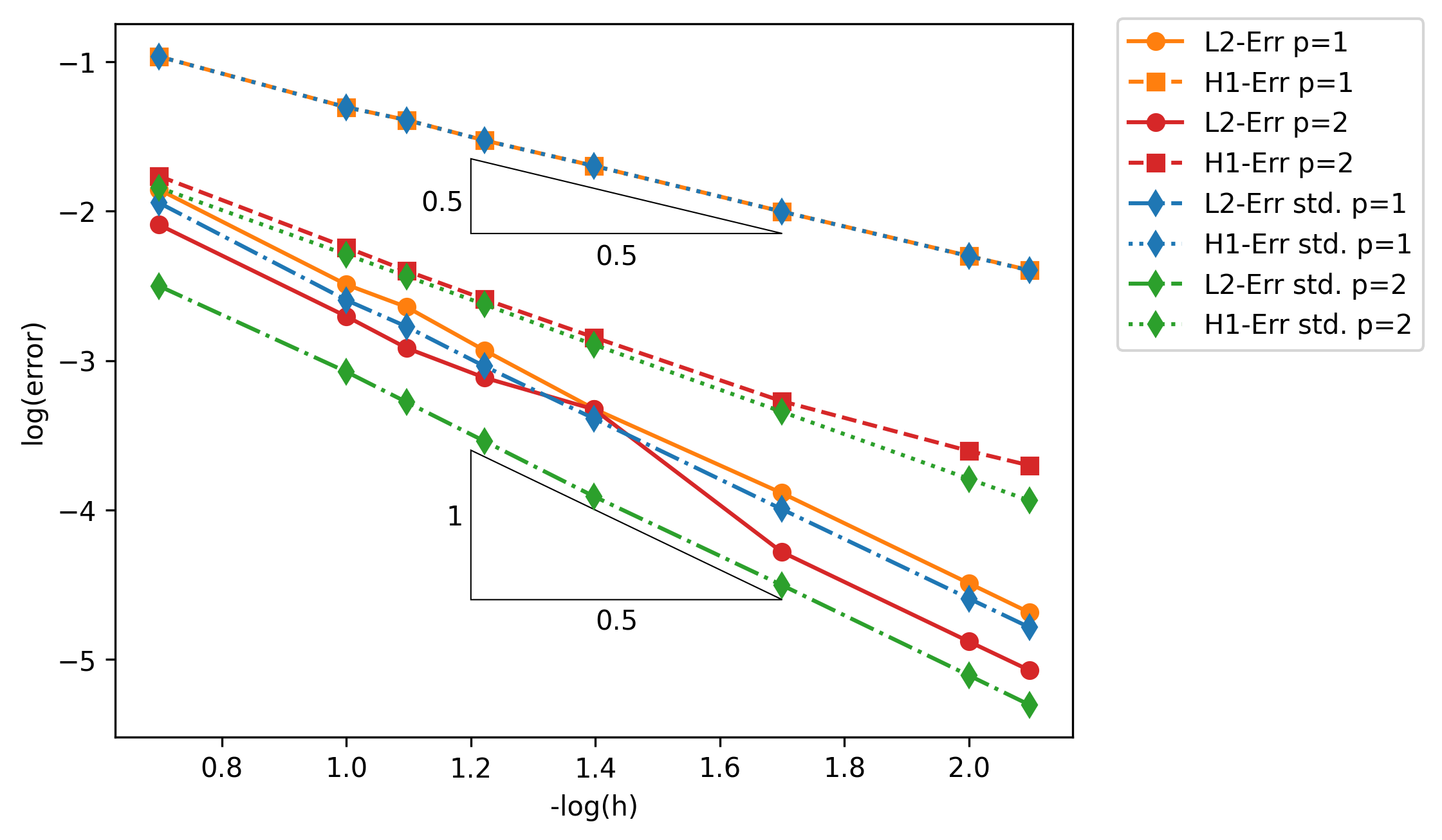}
		\caption{Error rates.}
		\label{fig:circ_usual}
     \end{subfigure}
     \hfill
     \begin{subfigure}[b]{0.45\textwidth}
        \centering
        \includegraphics[trim={0 0 0 0},clip,scale=0.45]{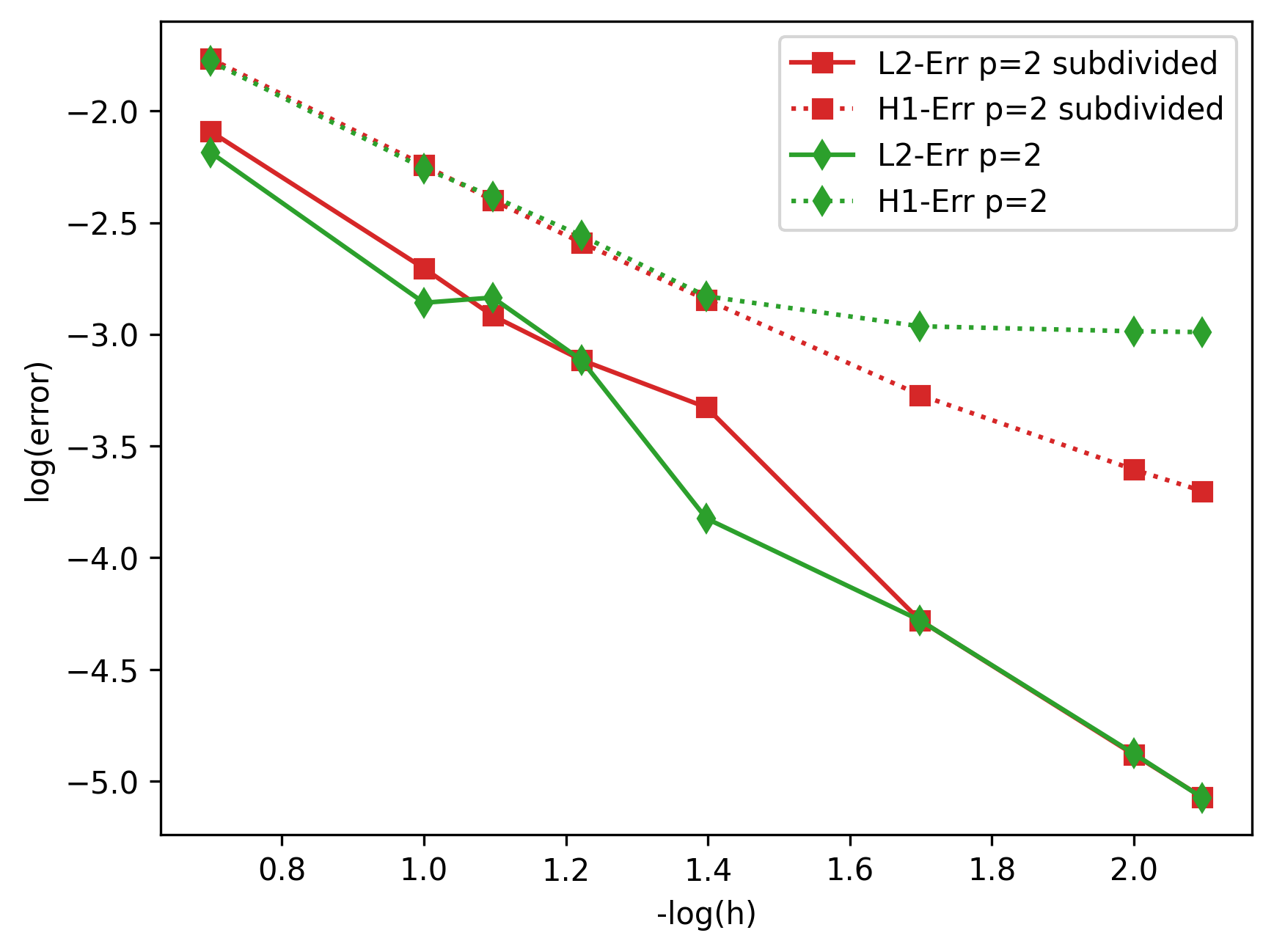}  
		\caption{Influence of the quadrature order.}
		\label{fig:circ_quad}
     \end{subfigure}
    \caption{Convergence for a disc shaped domain.}
\end{figure}

\begin{figure}
     \centering
     \begin{subfigure}[b]{0.45\textwidth}
         \centering
         \includegraphics[trim={0 0 0 0},clip,scale=0.5]{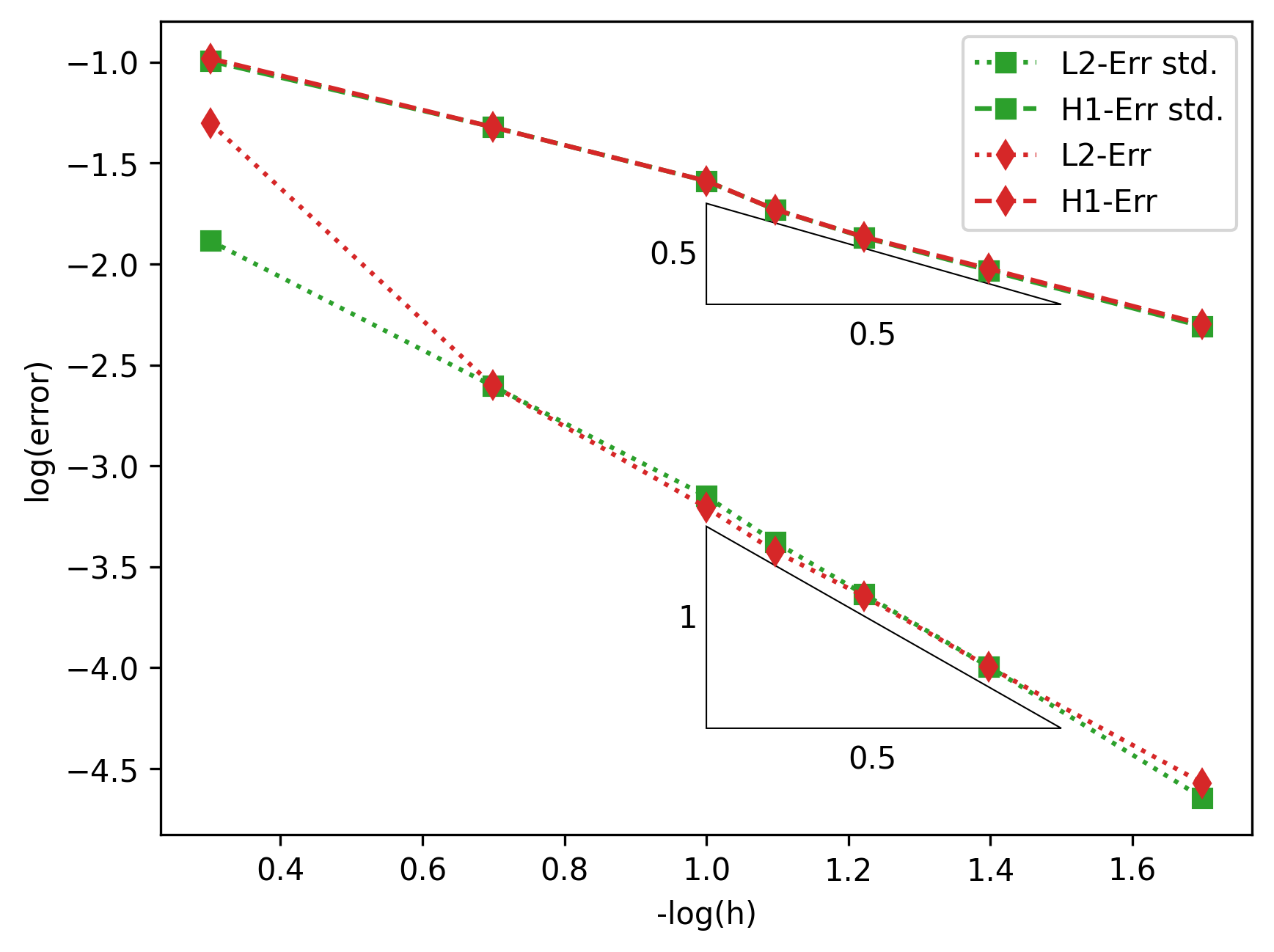}
		 \caption{Error rates.}
         \label{fig:triangle_usual}
     \end{subfigure}
     \hfill
     \begin{subfigure}[b]{0.45\textwidth}
         \centering
         \includegraphics[trim={5.5cm 0cm 5.5cm 1cm},clip,scale=0.25]{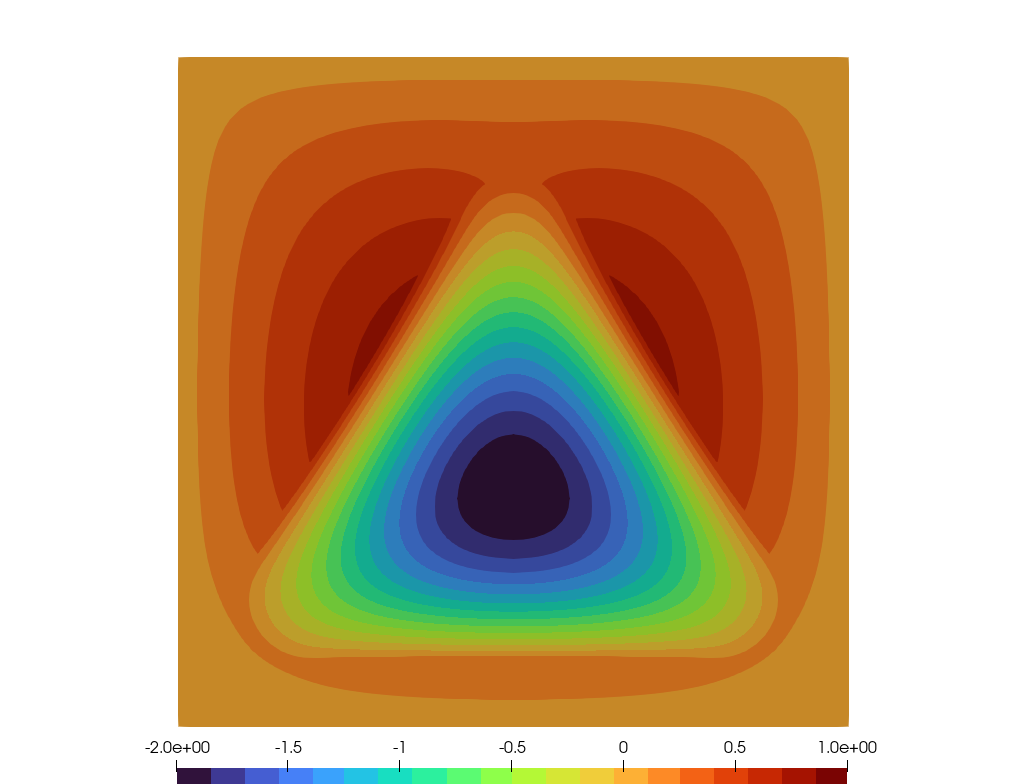} 
         \caption{Computed solution for \(h=0.04\) and contrast \(k_{-}=1.1\)}
         \label{fig:sol_triangle}
     \end{subfigure}
    \caption{Convergence for smooth triangle shaped domain.}
    \label{fig:three graphs}
\end{figure}


With this setup we now compute the solutions using \(H^{1}\)-conforming finite element spaces of order \(p=1\) and \(p=2\) for meshes with varying mesh sizes.
At the boundary and the interface we use curved elements with a quadratic geometry approximation.
In \Cref{fig:circ_quad} we compare for our new method the achieved errors with standard quadrature rules vs. quadrature rules with subdivision, where each element is divided into nine smaller similar elements.
We observe in \Cref{fig:circ_quad} that the adapted quadrature rule becomes necessary to achieve small \(H^{1}\)-errors.
Henceforth, this quadrature rule is used for all following calculations using the new method.

\Cref{fig:circ_usual} shows a log-log-plot of the relative errors in \(H^{1}\)- and in \(L^{2}\)-norm with respect to to the reference solution.
We compare our new method with a classical finite element method using the same meshes.
We observe that the errors are of comparable size and converge with the same rate.\\

After we have investigated the convergence for a simple domain, we will now move on to a more complicated domain to show and inspect the applicability of the method for more realistic configurations.
To this end we consider an equilateral triangle with rounded corners inside a square.
Here \(\Omega\) is the square with corners \((0,0),(10,0),(10,10)\) and \((0,10)\).
Now we consider the equilateral triangle \(\calD\) inside this square with corners \((2,2),(8,2)\) and \((5,2+3\sqrt{3})\) and define \(\dint{\Omega}:=\{\p\colon\dist(\p,\calD)<1\}\).
This leads to a shape with a boundary that is comprised of three circular arcs with radius \(1\) connected by three straight lines. As usual we then set \(\dext{\Omega}:=\Omega\setminus\overline{\dint{\Omega}}\) and we define \(\sigma\) to be piecewise constant such that \(\sigma|_{\dint{\Omega}}=-1\) and \(\sigma|_{\dext{\Omega}}=10\).
This enables us to choose \(\delta=0.5\). Finally, we choose $f=1$ as right hand-side and compare our solutions to a reference solution that was computed using the usual finite element method on a finer mesh (\(h=0.008\)). 

The resulting errors for finite elements of order \(1\) are depicted in \Cref{fig:triangle_usual}. There we observe that the usual method and our method converge with approximately the same rates, which shows that the new method also performs well for more complicated interface geometries.
In \Cref{fig:sol_triangle} we plot a computed solution for a contrast that is much closer to \(1\),
which is obtained by using \(\delta=h=0.04\).
While we have no analytical reference solution to compare with, we can still see the absence of singularities and the expected symmetry.

\subsection{Example in three dimensions}
\begin{figure}
     \centering
     \begin{subfigure}[b]{0.45\textwidth}
         \centering
        \includegraphics[trim={0 0 0 0},clip,scale=0.5]{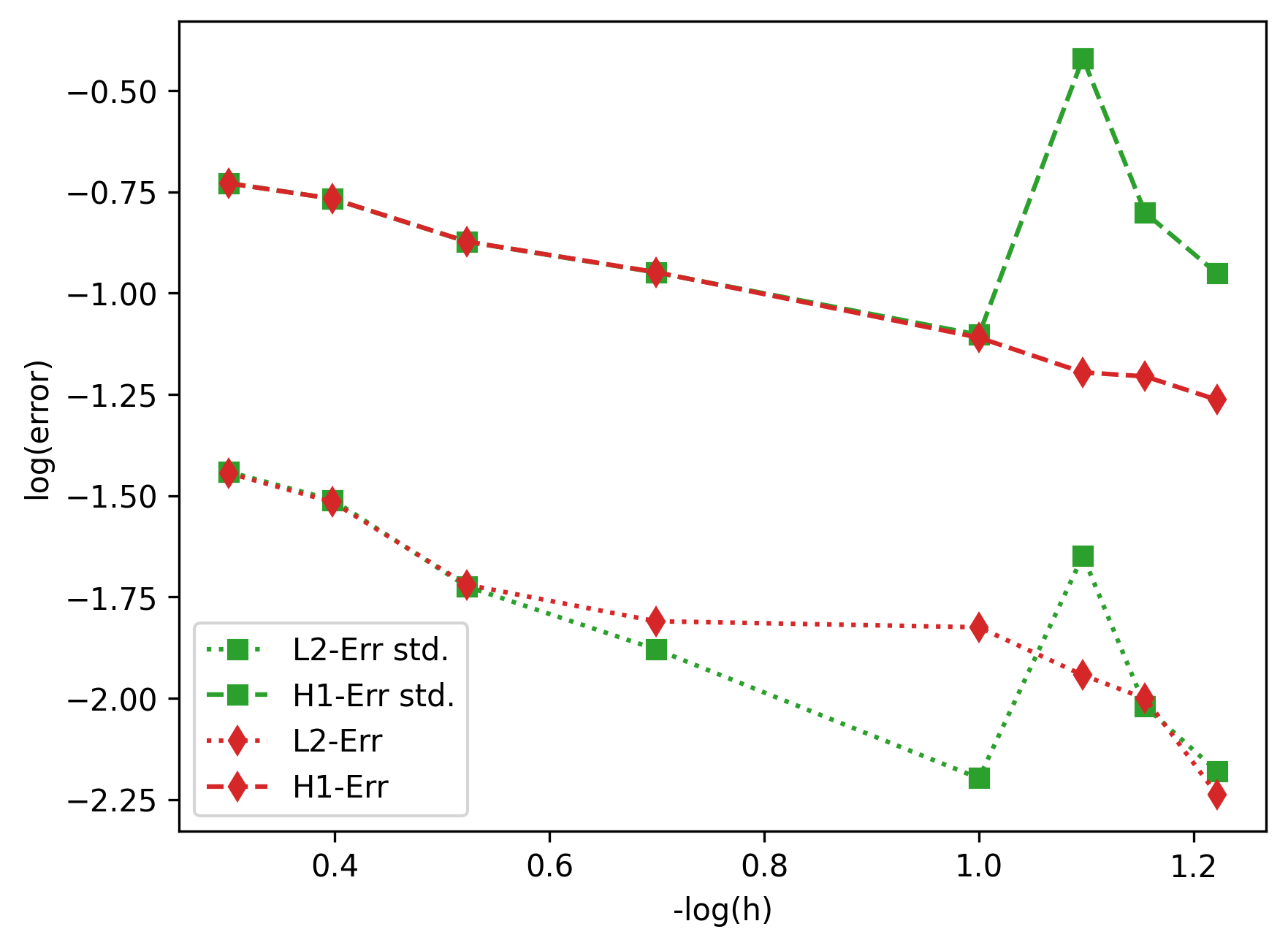}
		\caption{Errors rates.}
				\label{fig:ball_usual}
     \end{subfigure}
     \hfill
     \begin{subfigure}[b]{0.45\textwidth}
        \centering
        \includegraphics[clip,scale=0.23]{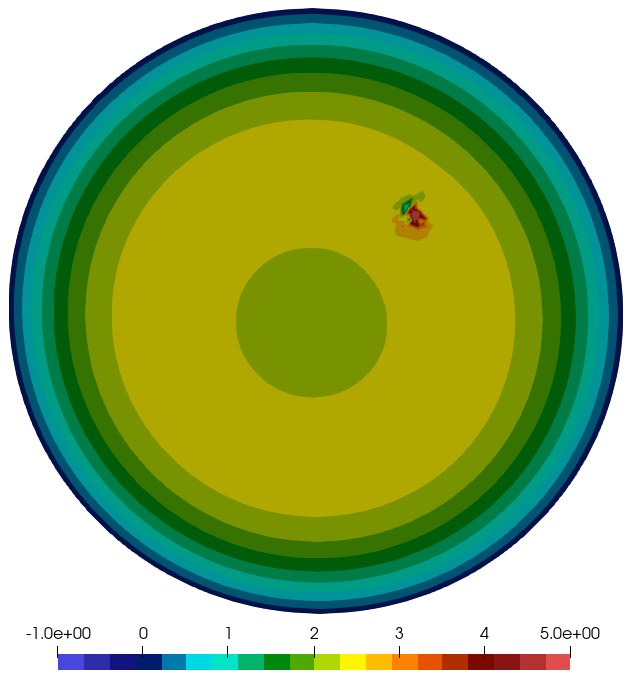}
		\caption{Cross-section of solution computed with the standard method \((h=0.08)\).}
				\label{fig:sol_wrong}
     \end{subfigure}
    \caption{Convergence for a ball shaped domain.}
\end{figure}

We also present an example for a three dimensional domain.
Again, to have an explicit reference solution we choose a domain that consist of a smaller ball inside a bigger one. This leads to the following domains
\begin{align*}
	\Omega&:=\{\p:|\p|< 4\},\qquad
	\dint{\Omega}:=\{\p:|\p|<2\},\qquad
	\dext{\Omega}:=\Omega\setminus\overline{\dint{\Omega}}.
\end{align*}
Thence we use the right hand-side
\[f(r):=\sigma(r)\frac{6r-9}{4}\]
corresponding to the solution
\[\trf_{\text{ref}}(r)=-\frac{1}{8}(r^{3}-3r^{2}-16)\]
in spherical coordinates, where \(\sigma|_{\dint{\Omega}}:=-1\) and \(\sigma|_{\dext{\Omega}}:=2\).
We then choose \(\delta=0.2\) and compute the relative errors which are depicted in \Cref{fig:ball_usual}.

\begin{figure}
     \centering
     \begin{subfigure}[b]{0.45\textwidth}
         \centering
        \includegraphics[trim={0 0 0 0},clip,scale=0.21]{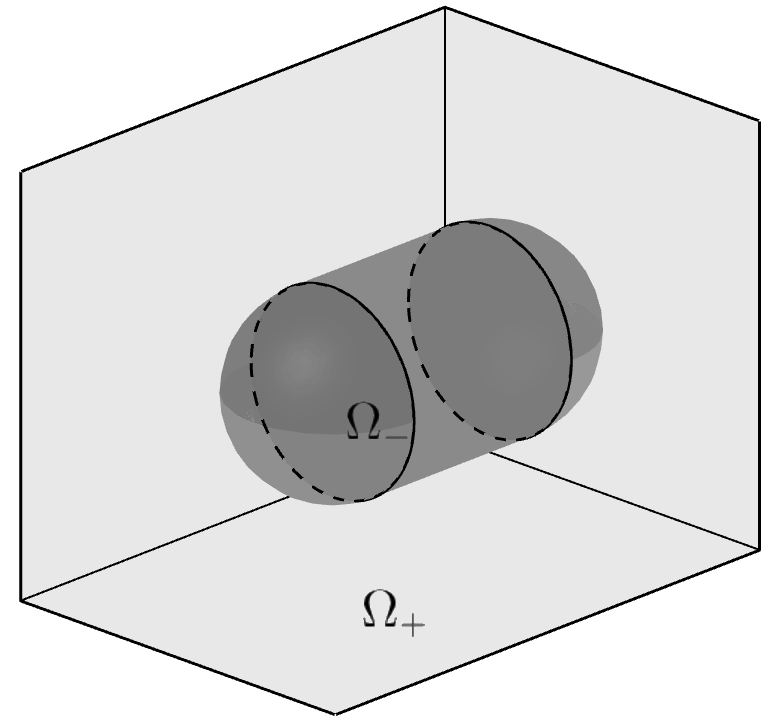} 
		\caption{Sketch of geometry.}
				\label{fig:pill_geo}
     \end{subfigure}
     \hfill
     \begin{subfigure}[b]{0.45\textwidth}
        \centering
        \includegraphics[trim={9.3cm 7cm 8cm 9cm},clip,scale=0.4]{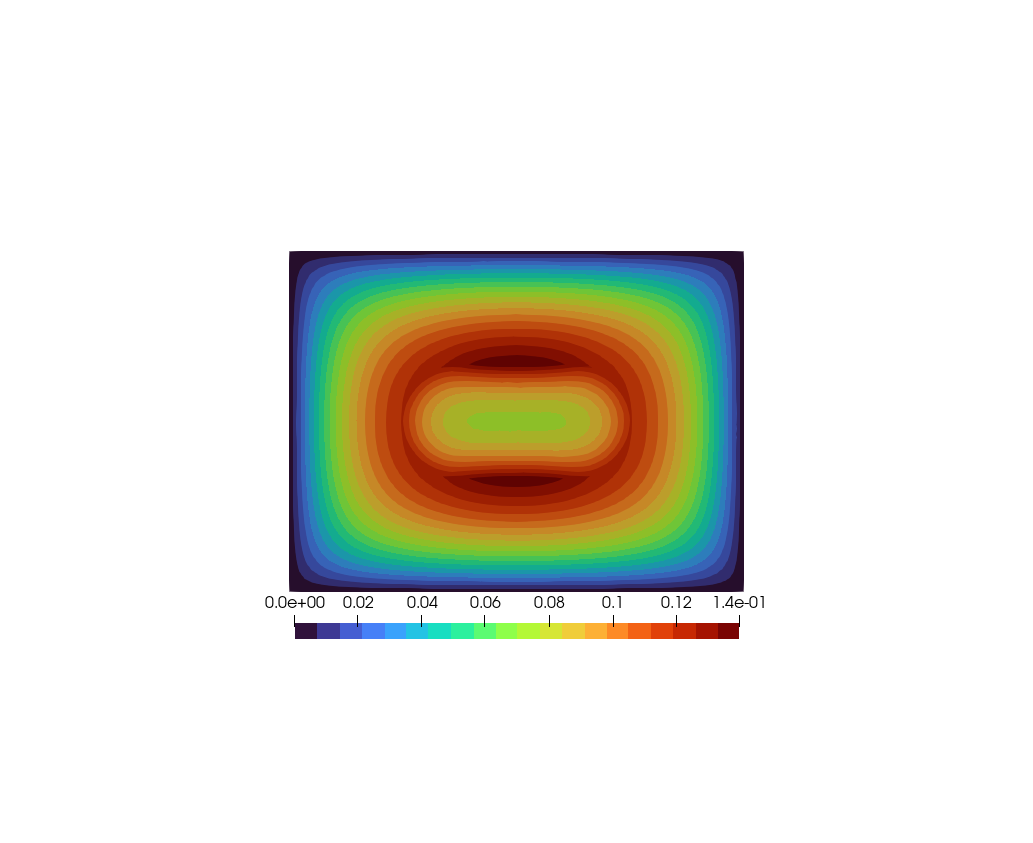}  
		\caption{Cross-section of the computed solution \((h=0.06)\).}
				\label{fig:sol_pill}
     \end{subfigure}
    \caption{Pill shaped inclusion.}
\end{figure}

We observe that the \(L^2\)- and \(H^{1}\)-errors for our method decay mostly as expected apart from a small plateau in the \(L^{2}\)-error, which is most likely caused by the use of anisotropic meshes for \(h>\delta\).
In contrast, we notice that the errors for the standard method do not converge. This phenomenon is caused by the appearance of local singularities near the interface depicted in \Cref{fig:sol_wrong}, which may occur regardless of the mesh size and have been observed previously, see, e.g., \cite{BonnetBDCarvalhoCiarlet:18}.
Finally, we also present an example for the more complicated pill shaped inclusion depicted in \Cref{fig:pill_geo}. A slice of the solution computed by our method can be seen in \Cref{fig:sol_pill}. We note that it has the expected symmetries and does not exhibit any singularities.

\subsection{Example of a dispersive eigenvalue problem}
\begin{figure}
     \centering
     \begin{subfigure}[b]{0.45\textwidth}
         \centering
        \includegraphics[trim={0 0 0 0},clip,scale=0.5]{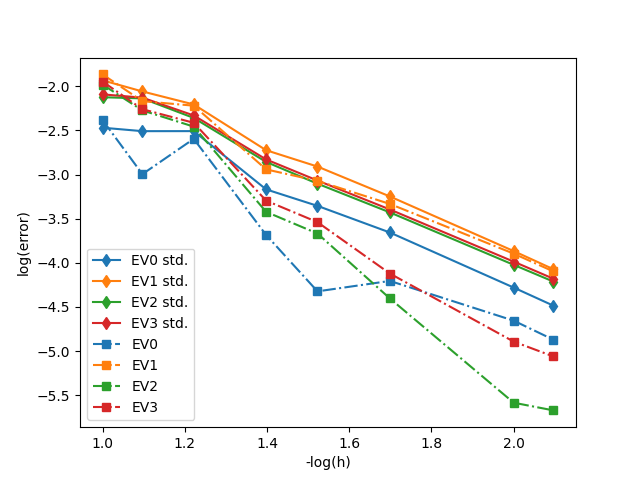}
		\caption{Errors for different eigenvalues.}
				\label{fig:conv_eig}
     \end{subfigure}
     \hfill
     \begin{subfigure}[b]{0.45\textwidth}
         \centering
        \includegraphics[clip,scale=4]{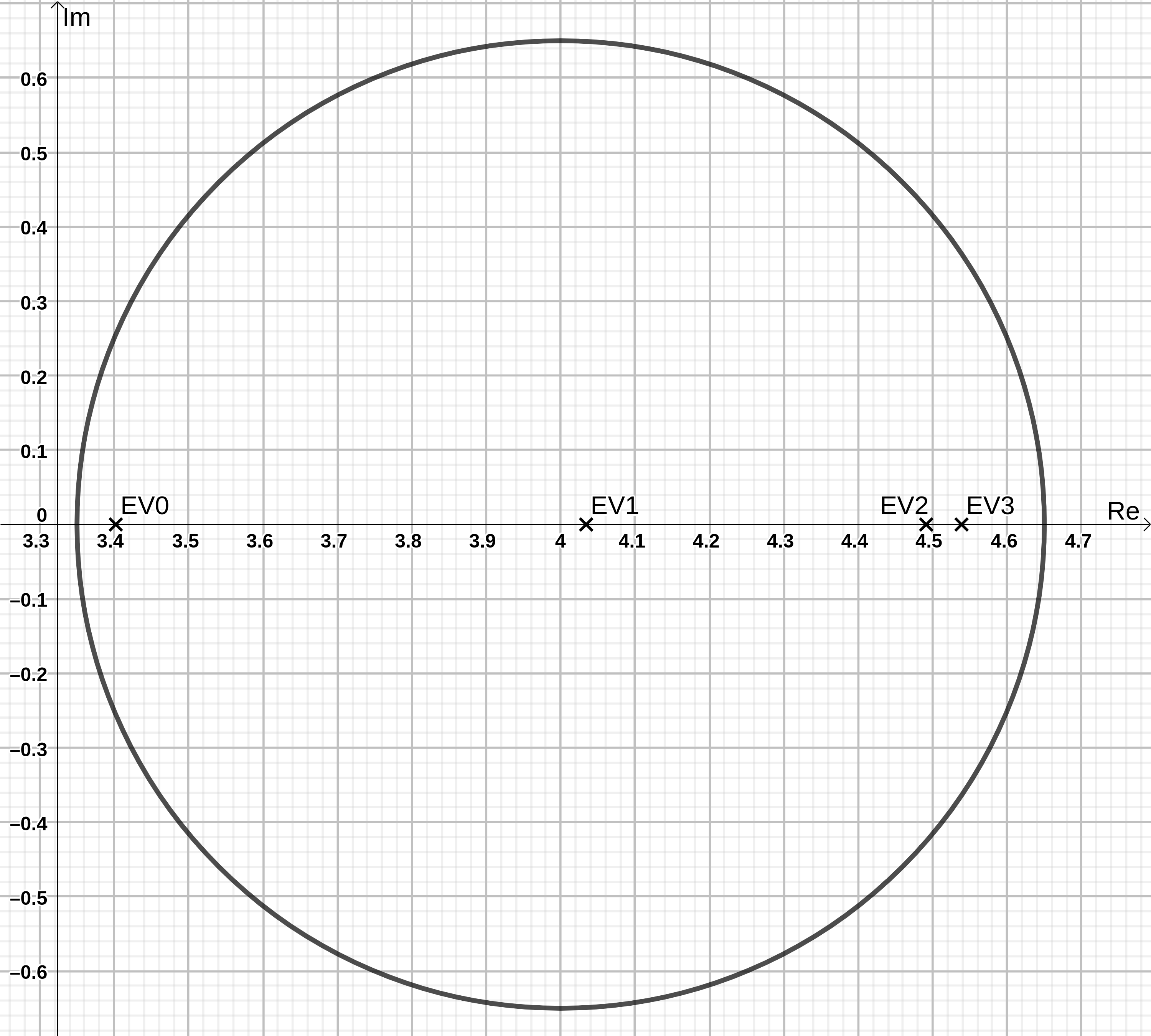} 
		\caption{Eigenvalues and the chosen contour.}
				\label{fig:cont}
     \end{subfigure}
    \caption{Dispersive eigenvalue problem.}
\end{figure}

We use the same disc shaped geometry, which we considered in the very first example.
Its symmetry allows us to compute the eigenvalues semi-analytically using Bessel functions and this enables us to obtain accurate reference solutions. We then consider the following eigenvalue problem:
\begin{align*}
	\text{Find } (\omega,\trf)\in \setC\times H_{0}^{1}(\Omega)\setminus\{0\} \text{ such that }
	-\div(\sigma(\omega)\nabla\trf){\color{blue}-}\omega^{2}\trf=0\\
\end{align*}
where
\begin{align*}
	\sigma(\omega)&:=\begin{cases}\frac{\omega^{2}}{\omega^{2}-200}& x\in \dint{\Omega}\\1& \text{else}\end{cases}.
\end{align*}
Our finite element method discretizes this problem into a holomorphic eigenvalue problem for a matrix, which is subsequently solved using the contour integral method proposed by Beyn \cite{Beyn:12}.
As the contour we choose a circle with radius $0.65$ centered at $4.0$. For comparison, the same method is also used for the discrete system obtained via a standard FEM. In \Cref{fig:conv_eig} we observe that both methods find the same 4 eigenvalues which are depicted in \Cref{fig:cont} and coincide with eigenvalues obtained by the semi-analytic method. We see that both methods reliably compute the eigenvalues with similar convergence rates.

\appendix
\section{Appendix}\label{sec:appendix}
In this appendix we present explicit bounds for the norms of the global reflection operators $R_\pm$, which are solely based on the curvature of the interface and \(\delta\).
The analysis consists of two main steps. First we prove the following lemma:
\begin{lemma}\label{lem:trafo_bound}
For \(R_{\pm}\) defined via \(\varphi\) as in \Cref{sec:global_reflection} we have
\[
	\|R_{\pm}\|\leq\sup_{\q\in \Phi^{-1}(\Sigma_{\mp})}|(D_{M(\q)}\Phi)(D_{\q}M)(D_{\q}\Phi)^{-1}|_{\calL(\setR^{d})}.
\]
\end{lemma}
\begin{proof}
We take \(w\in\calL(H_{0,\Gamma}^{1}(\Sigma_{\pm}))\) and then get by the chain rule
\begin{align*}
	||R_{\pm}w||_{H_{0,\Gamma}^{1}(\Sigma_{\mp})}^{2}&=\int_{\Sigma_{\mp}}{|\nabla(w\circ\Phi\circ M\circ\Phi^{-1})(\p)|^{2}}\dd \p\\
	&=\int_{\Sigma_{\mp}}{|[(\nabla w)(\Phi\circ M\circ\Phi^{-1}(\p))]^{T}(D_{M\circ\Phi^{-1}(\p)}\Phi)(D_{\Phi^{-1}(\p)}M)(D_{\p}\Phi^{-1})|^{2}}\dd \p\\
\end{align*}
Now we use the transformation formula for \(\q:=\Phi^{-1}(\p)\) and get
\begin{align*}
	||R_{\pm}w||_{H_{0,\Gamma}^{1}(\Sigma_{\mp})}^{2}&=\int_{\Phi^{-1}(\Sigma_{\mp})}{|[(\nabla w)(\Phi\circ M(\q))]^{T}(D_{M(\q)}\Phi)(D_{\q}M)(D_{\Phi(\q)}\Phi^{-1})|^{2}|\det D_{\q}\Phi|}\dd \q\\
	&=\int_{\Phi^{-1}(\Sigma_{\mp})}{|[(\nabla w)(\Phi\circ M(\q))]^{T}(D_{M(\q)}\Phi)(D_{\q}M)(D_{\q}\Phi)^{-1}|^{2}|\det D_{\q}\Phi|}\dd \q\\
&\leq	\int_{\Phi^{-1}(\Sigma_{\mp})}{|(\nabla w)(\Phi\circ M(\q))|^{2}||(D_{M(\q)}\Phi)(D_{\q}M)(D_{\q}\Phi)^{-1}||_{\calL(\setR^{d})}^{2}|\det D_{\q}\Phi|}\dd \q\\
\end{align*}
We define \(C_{\pm}:=\sup_{\q\in \Phi^{-1}(\Sigma_{\mp})}||(D_{M(\q)}\Phi)(D_{\q}M)(D_{\q}\Phi)^{-1}||_{\calL(\setR^{d})}\) and then apply the transformation formula again which yields
\begin{align*}
||R_{\pm}w||_{H_{0,\Gamma}^{1}(\Sigma_{\mp})}^{2}&\leq C_{\pm}^{2}\int_{\Phi^{-1}(\Sigma_{\mp})}{|(\nabla w)(\Phi\circ M(\q))|^{2}|\det D_{\q}\Phi|}\dd \q\\
&=C_{\pm}^{2}\int_{M\circ\Phi^{-1}(\Sigma_{\mp})}{|(\nabla w)(\Phi(\q'))|^{2}|\det D_{\q'}\Phi|}\dd \q'\\
&=C_{\pm}^{2}\int_{\Sigma_{\pm}}{|(\nabla w)(\p')|^{2}}\dd \p'=C_{\pm}^{2}||w||_{H_{0,\Gamma}^{1}(\Sigma_{\pm})}^{2}
\end{align*}
\end{proof}
As a second step we parameterize the interface and calculate the spectral norm of the matrix given by the lemma to get an explicit bound. This part is split into different sections based on the dimension.

\subsection{Bounds in two dimensions}\label{sec:bounds2d}
In case of two dimensions, the interface is a curve that is \(C^{1}\) and piece-wise \(C^{2}\).
Hence we can define the signed curvature \(\kappa\colon\Gamma\rightarrow\setR\) corresponding to the normal vector. 
We then obtain the following bounds for \(R_{\pm}\).
\begin{theorem}\label{the:bound2d}
In two dimensions the reflection operators are bounded in norm by
\begin{align*}
	\|\dint{R}\|&\leq\max\left(1,\left|\frac{1-\delta\inf_{\p\in\Gamma}\kappa(\p)}{1+\delta\inf_{\p\in\Gamma}\kappa(\p)}\right|\right),\qquad
	\|\dext{R}\|\leq\max\left(1,\left|\frac{1+\delta\sup_{\p\in\Gamma}\kappa(\p)}{1-\delta\sup_{\p\in\Gamma}\kappa(\p)}\right|\right).
\end{align*}
\end{theorem}
\begin{proof}
	We take \(\p\in \Gamma\) fixed and choose a coordinate system such that \(n(\p)=(0,1)^\top\).
Now we can parameterize \(\Gamma\) around \(\p\) by \(\gamma:[-\alpha,\alpha]\rightarrow\Gamma\) such that
\begin{align*}
	\gamma(0)=\p,\quad\gamma'(0)=(1,0)^\top. 
\end{align*}
Next we use the identity \(\frac{\dd}{\dd s}[n(\gamma(s))]=\kappa(\gamma(s))\gamma'(s)\) \cite[Chapter~1.5]{DoCarmo:76} to calculate \(D\Phi\) at \((\p,t_{0})\):
\begin{align*}
D_{(\p,t_{0})}\Phi&=\left(\begin{array}{c|c}\frac{\dd}{\dd s}[\gamma(s)+tn(\gamma(s))]_{s=0,t=t_{0}}&\frac{\dd}{\dd t}[\gamma(s)+tn(\gamma(s))]_{s=0,t=t_{0}}\end{array}\right)\\
&=\left(\begin{array}{c|c}\gamma'(0)+t_{0}\frac{\dd}{\dd s}[n(\gamma(s))]_{s=0}&n(\gamma(0))\end{array}\right)\\
&=\left(\begin{matrix}1+t_{0}\kappa(\p)&0\\0&1\end{matrix}\right).
\end{align*}
Because this formula is independent of the parametrization \(\gamma\), we can now calculate
\begin{align*}
	(D_{M((\p,t))}\Phi)(D_{(\p,t)}M)(D_{(\p,t)}\Phi)^{-1}&=(D_{(\p,-t)}\Phi)(D_{(\p,t)}M)(D_{(\p,t)}\Phi)^{-1}\\
	&=\left(\begin{matrix}1-t\kappa(\p)&0\\0&1\end{matrix}\right)\left(\begin{matrix}1&0\\0&-1\end{matrix}\right)\left(\begin{matrix}1+t\kappa(\p)&0\\0&1\end{matrix}\right)^{-1}\\
	&=\left(\begin{matrix}\frac{1-t\kappa(\p)}{1+t\kappa(\p)}&0\\0&-1\end{matrix}\right)
\end{align*}
and get
\begin{align*}
	|(D_{M((\p,t))}\Phi)(D_{(\p,t)}M)(D_{(\p,t)}\Phi)^{-1}|_{\calL(\setR^{2})}=\max\left(\left|\frac{1-t\kappa(\p)}{1+t\kappa(\p)}\right|,1\right).
\end{align*}
Using \Cref{lem:trafo_bound}, this leads to the following bound:
\begin{align*}
\|\dint{R}\|&\leq\sup_{(\p,t)\in\Phi^{-1}(\dext{\Sigma})}|(D_{M((\p,t))}\Phi)(D_{(\p,t)}M)(D_{(\p,t)}\Phi)^{-1}|_{\calL(\setR^{d})}\\
&=\sup_{(\p,t)\in\Phi^{-1}(\dext{\Sigma})}\left|\frac{1-t\kappa(\p)}{1+t\kappa(\p)}\right|
=\sup_{\p\in\Gamma,t\in[0,\delta)}\left|\frac{1-t\kappa(\p)}{1+t\kappa(\p)}\right|.
\end{align*}
If \(\kappa(\p)\geq0\), we always have that \(\left|\frac{1-t\kappa(\p)}{1+t\kappa(\p)}\right|\) attains its maximum at \(t=0\) and the maximum is \(1\). If \(\kappa(\p)\leq0\), the expression is monotonically increasing in \(t\) and its maximum is attained at \(t=\delta\). With this, we can eliminate the dependence of the supremum on \(t\) and obtain that
\begin{align*}
	\|\dint{R}\|&\leq\max\left(1,\left|\frac{1-\delta\inf_{\p\in\Gamma}\kappa(\p)}{1+\delta\inf_{\p\in\Gamma}\kappa(\p)}\right|\right).
\end{align*}
The result for \(\dext{R}\) can be derived in the same way.
\end{proof}

\subsection{Bounds in three dimensions}\label{sec:bounds3d}
Similar to the two dimensional case, we will again derive a bound for the operator based on the curvature of the surface. For this, we consider the two principal curvatures \(\kappa_{1},\kappa_{2}\colon\Gamma\mapsto\setR\) corresponding to the sign convention in the previous definition of the normal vector.
In the three dimensional case the map \(n\) is known as the Gauss map and we will use its properties to prove the following theorem.
\begin{theorem}\label{the:bound3d}
In three dimensions the reflection operators are bounded in norm by
\begin{align*}
\|\dint{R}\|&\leq\max\left(1,\left|\frac{1-\delta\inf_{\p\in\Gamma}\kappa_{1}(\p)}{1+\delta\inf_{\p\in\Gamma}\kappa_{1}(\p)}\right|,\left|\frac{1-\delta\inf_{\p\in\Gamma}\kappa_{2}(\p)}{1+\delta\inf_{\p\in\Gamma}\kappa_{2}(\p)}\right|\right),\\
\|\dext{R}\|&\leq\max\left(1,\left|\frac{1+\delta\sup_{\p\in\Gamma}\kappa_{1}(\p)}{1-\delta\sup_{\p\in\Gamma}\kappa_{1}(\p)}\right|,\left|\frac{1+\delta\sup_{\p\in\Gamma}\kappa_{2}(\p)}{1-\delta\sup_{\p\in\Gamma}\kappa_{2}(\p)}\right|\right).
\end{align*}
\end{theorem}
\begin{proof}
The proof is very similar to the two dimensional case. We again consider a fixed point \(p\in \Gamma\) and choose a coordinate system and a parametrization \(\gamma:(-\alpha,\alpha)\times(-\beta,\beta)\rightarrow\Gamma\) such that
\begin{align*}
	\gamma(0,0)=p,\quad\frac{\dd}{\dd x}\gamma(0,0)=\left(\begin{matrix}1\\0\\0\end{matrix}\right),\quad
	\frac{\dd}{\dd y}\gamma(0,0)=\left(\begin{matrix}0\\1\\0\end{matrix}\right),\quad
	n(p)=\left(\begin{matrix}0\\0\\1\end{matrix}\right).
\end{align*}
As in the two dimensional case, we then calculate
\begin{align*}
	D_{(p,t_{0})}\Phi&=\left(\begin{array}{c|c|c}\frac{\dd}{\dd x}[\gamma(x,y)+tn(\gamma(x,y))]&\frac{\dd}{\dd y}[\gamma(x,y)+tn(\gamma(x,y))]&\frac{\dd}{\dd t}[\gamma(x,y)+tn(\gamma(x,y))]\end{array}\right)\\
	&=\left(\begin{array}{c|c|c}\frac{\dd}{\dd x}\gamma(x,y)&\frac{\dd}{\dd y}\gamma(x,y)&n(\gamma(x,y))\end{array}\right)+t_{0}\left(\begin{array}{c|c|c}\frac{\dd}{\dd x}n(\gamma(x,y))&\frac{\dd}{\dd y}n(\gamma(x,y))&0\end{array}\right)\\
	&=\left(\begin{matrix}1&0&0\\0&1&0\\0&0&1\end{matrix}\right)+t_{0}\left(\begin{array}{c|c|c}\frac{\dd}{\dd x}n(\gamma(x,y))&\frac{\dd}{\dd y}n(\gamma(x,y))&0\end{array}\right).
\end{align*}
Now, we use that the derivative of the Gauss map \(n\) at \(p\) is given by the Weingarten map \(w\colon T_{p}\Gamma\mapsto T_{p}\Gamma\) where \(T_{p}\Gamma\) is the tangent space of \(\Gamma\) at \(p\).
Next we use that in our chosen coordinate system \(w\) can be represented by a matrix \(W\). Then we can write the second summand in the previous equation using the Weingarten map as
\begin{align*}
D_{(p,t_{0})}\Phi&=\left(\begin{matrix}1&0&0\\0&1&0\\0&0&1\end{matrix}\right)+t_{0}\left(\begin{array}{c|c}W&0\\\hline0&0\end{array}\right)
=\left(\begin{array}{c|c}I+t_{0}W&0\\\hline0&1\end{array}\right).
\end{align*}
We can now calculate
\begin{align*}
	(D_{M((p,t))}\Phi)(D_{(p,t)}M)(D_{(p,t)}\Phi)^{-1}&=\left(\begin{array}{c|c}I-tW&0\\\hline0&1\end{array}\right)\left(\begin{array}{c|c}I&0\\\hline0&-1\end{array}\right)\left(\begin{array}{c|c}I+tW&0\\\hline0&1\end{array}\right)^{-1}\\
	&=\left(\begin{array}{c|c}(I-tW)(I+tW)^{-1}&0\\\hline0&-1\end{array}\right).
\end{align*}
The Weingarten map is diagonalizable and its eigenvalues are the two principal curvatures of \(\Gamma\) at \(p\) \cite[Chapter~3.2]{DoCarmo:76}. We can therefore write it as
\(W=:SD_{w}S^{-1}\) with \[D_{w}:=\left(\begin{array}{cc}\kappa_{1}(p)&0\\0&\kappa_{2}(p)\end{array}\right)\]
and \(S\) being orthonormal. Inserting this into the equation above leads to
\begin{align*}
	(D_{M((p,t))}\Phi)(D_{(p,t)}M)(D_{(p,t)}\Phi)^{-1}&=\left(\begin{array}{c|c}(I-tSD_{w}S^{-1})(I+tSD_{w}S^{-1})^{-1}&0\\\hline0&-1\end{array}\right)\\
	&=\left(\begin{array}{c|c}S(I-tD_{w})S^{-1}S(I+tD_{w})^{-1}S^{-1}&0\\\hline0&-1\end{array}\right)\\
	&=\left(\begin{array}{c|c}S(I-tD_{w})(I+tD_{w})^{-1}S^{-1}&0\\\hline0&-1\end{array}\right)\\
	&=\tilde{S}\left(\begin{array}{c|c}(I-tD_{w})(I+tD_{w})^{-1}&0\\\hline0&-1\end{array}\right)\tilde{S}^{-1}\\
	&=\tilde{S}\left(\begin{matrix}\frac{1-t\kappa_{1}(p)}{1+t\kappa_{1}(p)}&0&0\\0&\frac{1-t\kappa_{2}(p)}{1+t\kappa_{2}(p)}&0\\0&0&-1\end{matrix}\right)\tilde{S}^{-1}
\end{align*}
with \[\tilde{S}:=\left(\begin{array}{c|c}S&0\\\hline0&1\end{array}\right).\]
This implies that
\[
	\|(D_{M((p,t))}\Phi)(D_{(p,t)}M)(D_{(p,t)}\Phi)^{-1}\|_{\calL(\setR^{3})}=\max\left(\left|\frac{1-t\kappa_{1}(p)}{1+t\kappa_{1}(p)}\right|,\left|\frac{1-t\kappa_{2}(p)}{1+t\kappa_{2}(p)}\right|,1\right).
\]
 Now the same calculations as in the two dimensional case lead to the claimed bounds.
\end{proof}

\subsection{Bounds for special geometries}
Finally, we apply the previous results to specific geometries which are used in our numerical experiments.
\begin{corollary}[Line or Plane]\label{cor:straight}
If the interface is a plane, the bounds are given by
 \(|R_{\pm}|\leq1.\)
\end{corollary}
\begin{corollary}[Part of a circle, sphere or cylinder]\label{cor:curved}
If the interface is a part of a sphere or cylinder with radius \(r\) and \(\delta<r\), we have to distinguish if the vector pointing inwards points towards \(\dint{\Omega}\) or \(\dext{\Omega}\).  Then the following bounds hold:
\begin{align*}
	\text{Towards }&\dint{\Omega}:&\|\dint{R}\|&\leq1&\text{ and }&&\|\dext{R}\|&\leq\frac{r+\delta}{r-\delta},\\
	\text{Towards }&\dext{\Omega}:&\|\dint{R}\|&\leq\frac{r+\delta}{r-\delta}&\text{ and }&&\|\dext{R}\|&\leq1.
\end{align*}
\end{corollary}
\begin{proof}
For the simple geometries considered in the corollaries above, the principal curvatures are constant and either \(\pm \frac{1}{r}\) or \(0\). By inserting these values into the bound from \Cref{the:bound3d} the given bounds then follow immediately.
\end{proof}
Note that it is possible to combine these different parts to achieve interfaces which are rounded polygons or polyhedra.

\bibliographystyle{abbrv}
\bibliography{short_biblio}
\end{document}